\documentclass[final,onefignum,onetabnum]{mousaviart}
\usepackage{enumerate}
\usepackage{mathrsfs}
\usepackage{stackrel}
\usepackage{bigints}
\usepackage{enumitem}
\newtheorem{assumption}{Assumption}
\usepackage{braket,amsfonts}

\usepackage{array}

\usepackage[caption=false]{subfig}
\usepackage{pgfplots}
\newsiamremark{remark}{Remark}
\usepackage{algorithmic}

\usepackage{graphicx,epstopdf}
\Crefname{ALC@unique}{Line}{Lines}

\usepackage{amsopn}

\usepackage{xspace}
\usepackage{bold-extra}
\usepackage[most]{tcolorbox}

\colorlet{texcscolor}{blue!50!black}
\colorlet{texemcolor}{red!70!black}
\colorlet{texpreamble}{red!70!black}
\colorlet{codebackground}{black!25!white!25}

\lstdefinestyle{siamlatex}{%
	style=tcblatex,
	texcsstyle=*\color{texcscolor},
	texcsstyle=[2]\color{texemcolor},
	keywordstyle=[2]\color{texemcolor},
	moretexcs={cref,Cref,maketitle,mathcal,text,headers,email,url},
}

\tcbset{%
	colframe=black!75!white!75,
	coltitle=white,
	colback=codebackground, % bottom/left side
	colbacklower=white, % top/right side
	fonttitle=\bfseries,
	arc=0pt,outer arc=0pt,
	top=1pt,bottom=1pt,left=1mm,right=1mm,middle=1mm,boxsep=1mm,
	leftrule=0.3mm,rightrule=0.3mm,toprule=0.3mm,bottomrule=0.3mm,
	listing options={style=siamlatex}
}

\newtcblisting[use counter=example]{example}[2][]{%
	title={Example~\thetcbcounter: #2},#1}

\newtcbinputlisting[use counter=example]{\examplefile}[3][]{%
	title={Example~\thetcbcounter: #2},listing file={#3},#1}

\DeclareTotalTCBox{\code}{ v O{} }
{ %fontupper=\ttfamily\color{texemcolor},
	fontupper=\ttfamily\color{black},
	nobeforeafter,
	tcbox raise base,
	colback=codebackground,colframe=white,
	top=0pt,bottom=0pt,left=0mm,right=0mm,
	leftrule=0pt,rightrule=0pt,toprule=0mm,bottomrule=0mm,
	boxsep=0.5mm,
	#2}{#1}

\patchcmd\newpage{\vfil}{}{}{}
\begin{tcbverbatimwrite}{tmp_\jobname_header.tex}
	\title{Mixed finite element approximation for\\ non-divergence form elliptic equations\\ with random input data}
			
	\author{\vspace*{5mm}Amireh Mousavi}
		
	\headers{Mixed FEM for non-divergence PDEs with random data}{A. Mousavi}
\end{tcbverbatimwrite}
	\title{Mixed finite element approximation for\\ non-divergence form elliptic equations\\ with random input data}
			
	\author{\vspace*{5mm}Amireh Mousavi}
		
	\headers{Mixed FEM for non-divergence PDEs with random data}{A. Mousavi}

\ifpdf
\hypersetup{ pdftitle={Mixed FEM for non-divergence form PDEs with uncertainity} }
\fi

\usepackage{lettrine}

\makeatletter

\def\GetFirst#1#2\relax{#1}
\def\GetRest#1#2\relax{#2}

\renewcommand{\section}[1]{%
	\refstepcounter{section}%
	\addcontentsline{toc}{section}{\protect\numberline{\thesection}#1}%
	\vspace{2ex}%
	\edef\temp{#1}%
	\edef\FirstLetter{\expandafter\GetFirst\temp\relax}%
	\edef\RestOfTitle{\expandafter\GetRest\temp\relax}%
	\begin{center}
		{\footnotesize
			\normalfont
			{\fontsize{10}{10}\selectfont \thesection.}\;
			{\fontsize{10}{10}\selectfont \MakeUppercase{\FirstLetter}}%
			\MakeUppercase{\RestOfTitle}%
		}
	\end{center}%
	\vspace{1.5ex}%
}
\makeatother
\makeatletter
\renewcommand\subsection[1]{%
	\refstepcounter{subsection}%
	\addcontentsline{toc}{subsection}{\protect\numberline{\thesection.\arabic{subsection}}#1}%
	\vspace{0.3ex plus 0.5ex minus 0.2ex}%
	\noindent
	{\normalfont\normalsize
		\thesection.\arabic{subsection}.\kern0.4em % فاصله دقیق بعد از شماره
		\bfseries #1.}%
}
\makeatother
\usepackage{changepage} % برای تغییر عرض متن
\usepackage{xcolor}

\newcommand{\FancyAbstractTitle}[1]{%
	{\fontsize{9pt}{1pt}\selectfont #1} % حرف اول بزرگتر
}
\begin{document}
	\maketitle
	\renewcommand{\footnoterule}{%
		\kern -3pt%-3pt
		\hrule width 0.2\linewidth height 0.2pt
		\kern 4pt
	}
	\begingroup
	\renewcommand\thefootnote{}%
	\footnote{This work was supported by the European Research Council (ERC Project DAFNE, grant agreement No. 891734).
	\\
		The author thanks Professor Dietmar Gallistl for valuable discussions and constructive comments that helped improve this work.}%
	\addtocounter{footnote}{-1}%
	\endgroup
	
	%% ------------------------------------------------------------------
	%% ABSTRACT
	%% ------------------------------------------------------------------
	\begin{tcbverbatimwrite}{tmp_\jobname_abstract.tex}
		%\begin{abstract}
		\begin{center}
			\begin{adjustwidth}{1cm}{1cm} % عرض کمتر از متن اصلی
				\noindent
				{\FancyAbstractTitle{A}{\hspace*{-.8mm}}\fontsize{6pt}{1pt}\selectfont{BSTRACT.}}\quad
				\footnotesize
				%{\fontsize{8pt}{1pt}\selectfont
			We consider an elliptic partial differential equation in non-divergence form with a random diffusion matrix and random forcing term.
			To address this, we propose a mixed-type continuous finite element discretization in the physical domain, combined with a collocation discretization in the stochastic domain.
			For the mixed formulation, we first introduce a stochastic cost functional at the continuous level. This formulation is then enhanced to incorporate the vanishing tangential trace constraint directly into a mesh-dependent cost functional, rather than enforcing it in the solution’s function space. In this context, we define a mesh-dependent norm and provide an error analysis based on this norm.
			We employ the collocation method by collocating the stochastic equation at the zeros of suitable tensor product orthogonal polynomials. This approach leads to a system of uncoupled deterministic problems, simplifying computation.
			Furthermore, we establish an a poriori error bound for the fully discrete approximation, detailing the convergence rates with respect to the discretization parameters. Finally, numerical results are presented to confirm and validate the theoretical findings.
		%}
		\end{adjustwidth}
	\end{center}
		%\end{abstract}
		
%		\begin{keywords}
%			elliptic PDEs with random data, non-divergence form, finite element method, stochastic collocation method, a priori error estimate
%		\end{keywords}
%		
%		\begin{MSCcodes}
%			65C20, 65N30, 65N35, 65N15
%		\end{MSCcodes}
	\end{tcbverbatimwrite}
			%\begin{abstract}
		\begin{center}
			\begin{adjustwidth}{1cm}{1cm} % عرض کمتر از متن اصلی
				\noindent
				{\FancyAbstractTitle{A}{\hspace*{-.8mm}}\fontsize{6pt}{1pt}\selectfont{BSTRACT.}}\quad
				\footnotesize
				%{\fontsize{8pt}{1pt}\selectfont
			We consider an elliptic partial differential equation in non-divergence form with a random diffusion matrix and random forcing term.
			To address this, we propose a mixed-type continuous finite element discretization in the physical domain, combined with a collocation discretization in the stochastic domain.
			For the mixed formulation, we first introduce a stochastic cost functional at the continuous level. This formulation is then enhanced to incorporate the vanishing tangential trace constraint directly into a mesh-dependent cost functional, rather than enforcing it in the solution’s function space. In this context, we define a mesh-dependent norm and provide an error analysis based on this norm.
			We employ the collocation method by collocating the stochastic equation at the zeros of suitable tensor product orthogonal polynomials. This approach leads to a system of uncoupled deterministic problems, simplifying computation.
			Furthermore, we establish an a poriori error bound for the fully discrete approximation, detailing the convergence rates with respect to the discretization parameters. Finally, numerical results are presented to confirm and validate the theoretical findings.
		%}
		\end{adjustwidth}
	\end{center}
		%\end{abstract}
		
%		\begin{keywords}
%			elliptic PDEs with random data, non-divergence form, finite element method, stochastic collocation method, a priori error estimate
%		\end{keywords}
%		
%		\begin{MSCcodes}
%			65C20, 65N30, 65N35, 65N15
%		\end{MSCcodes}

	%% ------------------------------------------------------------------
	%% END HEADER
	%% ------------------------------------------------------------------
	%
	\vspace*{6mm}
	\section{Introduction}
	\label{sec:intro}
	
	This work addresses the numerical analysis of elliptic stochastic partial differential equations of the form  
	\begin{equation}
		\label{eq:introduction}
		\mathcal{L}(\omega, \cdot) u(\omega, \cdot) = f(\omega, \cdot) \quad \mathrm{in}~ \mathcal{D},
	\end{equation}
	for almost surely (a.s.) $\omega \in \Omega$. Here, the linear elliptic operator $\mathcal{L}(\omega, \cdot)$ is expressed in the general form $\boldsymbol{A}(\omega, \cdot) : \mathrm{D}^2$, where $\boldsymbol{A}(\omega, \boldsymbol{x})$ represents a stochastic diffusion matrix and $\mathrm{D}^2$ is the Hessian operator with respect to the spatial variable $\boldsymbol{x}$. Both the diffusion matrix $\boldsymbol{A}(\omega, \boldsymbol{x})$ and the forcing term $f(\omega, \boldsymbol{x})$ are defined over the sample space $\Omega$ and the physical domain $\mathcal{D}$, making them stochastic functions. For simplicity and clarity, this study focuses on operators without lower-order terms, though the methodology and results can be readily extended to include such terms.
	
	Input uncertainties may arise from incomplete knowledge, which could theoretically be addressed with better measurements or devices. However, such remedies are often too costly or impractical. Additionally, some uncertainties come from inherent system variability that cannot be reduced through further experimentation or improved measurement tools.
	These sources of uncertainty frequently lead to the appearance of equations like (\ref{eq:introduction}), either directly or as a result of the linearization of nonlinear equations in various applications. For example, the linearization of the fully nonlinear Hamilton-Jacobi-Bellman equation produces a sequence of linear equations in non-divergence form \cite{Smears-Suli2, Gallistl-Süli, Lakkis-Mousavi2}. When the control maps are uncertain, this linearization results in equations of the type (\ref{eq:introduction}). A practical application of this can be seen in game theory with diffusion control \cite{Ankirchner-Kazi-Wendt-Zhou}, where the strategies of players are modeled as random fields rather than deterministic functions.
	
	The equation's domain consists of two distinct components: the physical domain $\mathcal{D}$ and the sample space $\Omega$. Consequently, it is crucial to employ appropriate discretization techniques for each component. In this work, we combine the finite element method, a well-established numerical approach for solving deterministic partial differential equations, with a suitable approximation strategy designed for the stochastic variables. This hybrid approach provides a robust and efficient framework for approximating solutions to the stochastic partial differential equation under consideration.
	
	When the operator $\mathcal{L}(\omega, \cdot)$ is in divergence form, such problems have been extensively studied in the literature. Various numerical methods have been developed and analyzed, including Monte Carlo type methods \cite{Grham&etal, Dick&etal1, Dick&etal2}, stochastic collocation methods \cite{Babuska-Nobile-Tempone, Nobile&Tempone&Webster, Xiu-Hesthaven}, and stochastic Galerkin methods, which are inherently intrusive, as discussed in \cite{Babuska-et-al, Frauenfelder-Schwab, Matthies-Keese}. Additionally, a stochastic Galerkin mixed formulation has been investigated in \cite{Ern&Pow&Silv&Ull}.
	
	In this paper, we focus on the case where the operator $\mathcal{L}(\omega, \cdot)$ is in non-divergence form with respect to the spatial variable $\boldsymbol{x}$. The lack of a natural variational structure precludes a straightforward use of weak solutions in $H^1(\mathcal{D})$ over the physical domain. However, under the Cordes condition assumed in Assumption~\ref{assum:A-f}~\textit{\ref{assum:matrix-A}\ref{assum:Cordes-condition}}, the well-posedness of a strong solution $u(\omega, \cdot) \in H^2(\mathcal{D})$ for a.s. $\omega \in \Omega$ is guaranteed. To leverage the computational advantages of working in the $H^1$ space when addressing the numerical approximation in the physical domain, we propose a least-squares approach combined with a gradient recovery method introduced in \cite{Gallistl} and further studied in \cite{Lakkis-Mousavi}. Specifically, we adopt a least-squares Galerkin gradient recovery method, which can be interpreted as a mixed finite element technique.
	
	The least-squares approach provides a significant advantage by replacing the constraints typically required for the well-posedness of a problem with additional terms in the quadratic functional. Specifically, we define a cost functional $E$ corresponding to the problem under consideration at the continuous level. This flexibility is especially useful when constructing finite element approximations, as it eliminates the need to enforce conditions such as curl-free constraints or vanishing tangential traces, which can be challenging to satisfy exactly.
	Handling vanishing tangential traces poses a significant challenge in establishing the coercivity of the associated bilinear form. To overcome this, we introduce a mesh-dependent cost functional $E_h$ along with a corresponding mesh-dependent norm. By using an inverse trace inequality in finite element spaces, we effectively resolve this issue within the discrete finite element framework.
	
	For problems where the inputs and, consequently, the solutions are sufficiently smooth functions of random variables, the stochastic collocation method (SCM) serves as a powerful alternative to traditional Monte Carlo sampling techniques. This method combines the strengths of both spectral and sampling approaches. On one hand, it preserves the accuracy of the spectral Galerkin method, achieving exponential convergence when the coefficient $\boldsymbol{A}$ and the forcing term $f$ are smooth with respect to the random variables. On the other hand, it offers the flexibility of sampling methods by decoupling the system of linear equations with respect to the random variables.
	
	In this work, we propose a collocation method, as presented in \cite{Babuska-Nobile-Tempone}, which involves collocating the problem at the zeros of tensor product orthogonal polynomials with respect to an auxiliary joint probability density $\hat{\rho}$, corresponding to independent random variables. However, the tensor product SCM is affected by the so-called curse of dimensionality, where the dimension of the approximating space grows exponentially with the number of random variables. A potential solution to this challenge is the sparse grid method \cite{Nobile&Tempone&Webster} and its anisotropic version \cite{Nobile&Tempone&Webster2, Beck-Nobile-Tamellini-Tempone}. In this context, the use of tensor product spaces is particularly advantageous when the number of random variables is small. This paper focuses on the case where the probability space has low dimensionality, meaning the stochastic problem depends on a relatively small number of random variables. 
	
	For the stochastic approximation, we will consider the SCM \cite{Babuska-Nobile-Tempone}. For discretization in the physical space, we will introduce a mixed-type finite element method, evaluated with a mesh-dependent norm.
	
	The structure of the paper is as follows:
	In Section~\ref{sec:Problem-setting}, we introduce the problem, present the necessary background material and notations, and establish the well-posedness of the strong solution to the equation.
	In Section~\ref{sec:Mixed-formulation}, we introduce a stochastic quadratic functional based on a mixed formulation and then consider its corresponding minimization problem, where the minimizer coincides with the solution to the main problem.	
	In Section~\ref{sec:discretization}, we discuss the discretization techniques for both the spatial and stochastic domains. For the spatial domain, we employ the mixed continuous finite element method, introducing a mesh-dependent cost functional and establishing the well-posedness of the discretization within a corresponding mesh-dependent norm. For the stochastic domain, we adapt the collocation method introduced by \cite{Babuska-Nobile-Tempone} to achieve an efficient and accurate approximation.	
	In Section~\ref{sec:A-priori-error-convergence}, we impose additional assumptions on the data to ensure the regularity of the solution with respect to the random variables and present an a priori error bound for the fully discrete approximation, demonstrating the convergence rate of the discretization.
	In Section~\ref{sec:Numerical-Experiment}, we provide numerical results to validate the theoretical findings.
	
	\vspace*{2mm}
	\section{Problem setting}
	\label{sec:Problem-setting}
	
	Let $\mathcal{D}$ denote a bounded convex domain in $\mathbb{R}^d$, $d=2,3$ and let $(\Omega, \mathcal{F}, P)$ represent a complete probability space. Here $\Omega$ is the set of outcomes, $\mathcal{F}\subseteq 2^{\Omega}$ is the $\sigma$-algebra of events, and $P:\mathcal{F}\rightarrow [0,1]$ is a probability measure. We consider the following linear elliptic stochastic boundary value problem (SBVP): find a stochastic function $u :  \Omega \times \bar{\mathcal{D}} \rightarrow \mathbb{R} $, 
	such that, almost surely (a.s.), the following equation is satisfied
	\begin{equation}
		\label{eq:nondivergence}
		\boldsymbol{A}(\omega, \cdot): \mathrm{D}^2u(\omega, \cdot) = f(\omega, \cdot) ~~ \mathrm{in}~\mathcal{D},
		\quad \mathrm{and} ~~
		u(\omega, \cdot) =0 ~~ \mathrm{on}~ \partial \mathcal{D}.
	\end{equation}
	Here $\boldsymbol{A}:\Omega \times \mathcal{D}  \rightarrow \text{sym}(\mathbb{R}^{d \times d})$ and $f: \Omega \times \mathcal{D} \rightarrow \mathbb{R}$ are stochastic functions.
	Since the sources of stochasticity of the random fields $\boldsymbol{A}(\omega, \boldsymbol{x})$ and $f(\omega, \boldsymbol{x})$ are often unrelated to each other, we may consider two independent probability spaces $(\Omega_{A}, \mathcal{F}_{A}, P_{A})$ and $(\Omega_f, \mathcal{F}_f, P_f)$ respectively for them. The solution $u$ is defined on the product probability space $ (\Omega, \mathcal{F}, P) = \left( \Omega_{A} \times \Omega_f, \mathcal{F}_{A} \times \mathcal{F}_{f}, P_{A} \times P_{f}  \right) $ and $\omega = (\omega_{A}, \omega_{f}) \in \Omega$, in which $\omega_{A} \in \Omega_{A}$ and $\omega_{f} \in \Omega_{f}$. In fact, $\boldsymbol{A}$ and $f$ are essentially functions of $\omega_{A}$ and $\omega_f$ respectively.
	To guarantee the well-posedness of the strong solution, we impose the following assumptions.
	
	\begin{assumption}
		\label{assum:A-f}
		\textcolor{white}{.}
		\begin{enumerate}[label=\textnormal{(\alph*)}]
			\item 
			\label{assum:matrix-A}
			The matrix coefficient $\boldsymbol{A}$ satisfies the following:
			\begin{enumerate}[label=\textnormal{(\roman*)}]
				\item
				\label{assum:uniform-ellipticity}
				\textnormal{(uniform ellipticity condition)}.
				There exists $\lambda \in (0,1]$ such that
				\begin{multline}
					\label{con:ellipticity}
					P_A
					\Big(
					\omega_A \in \Omega_A : \lambda  \lvert\boldsymbol{\xi}\rvert^2 \leq \boldsymbol{\xi}^T \boldsymbol{A}( \omega_A, \boldsymbol{x}) \boldsymbol{\xi}\leq \lambda^{-1} 
					\lvert\boldsymbol{\xi}\rvert^2 ~~ 
					\forall \boldsymbol{\xi} \in \mathbb{R}^d 
					\\
					\text{ and for almost every (a.e.) } 
					\boldsymbol{x} \in \mathcal{D}
					\Big) =1.
				\end{multline}
				\item
				\label{assum:Cordes-condition}
				\textnormal{(Cordes condition)}. 
				There exists $\varepsilon \in (0,1]$ such that 
				\begin{equation}
					\label{con:cordes}
					P_A \left( 
					\omega_A \in \Omega_A : \dfrac{\lvert \boldsymbol{A}(\omega_A, \boldsymbol{x}) \rvert^2}{\left( \mathrm{tr} {\boldsymbol{A}}(\omega_A, \boldsymbol{x}) \right) ^2}
					\leq
					\dfrac{1}{d-1+\varepsilon} ~~ 
					\text{for a.e. } \boldsymbol{x} \in \mathcal{D}
					\right) = 1. 
				\end{equation}
			\end{enumerate}
			\item 
			The right-hand side random field $f \in L^2_{P_f}(\Omega_f, L^2(\mathcal{D}))$, which means
			\begin{equation}
				\label{assum:on-f-expection}
				\int_{\Omega_f}{\int_{\mathcal{D}} f^2(\omega_f, \boldsymbol{x})~d \boldsymbol{x}} d P_f(\omega_f) 
				< \infty.
			\end{equation}
			Equivalently,
			\begin{equation}
				{\label{assum:on-f}
					P_f\left(  \omega_f \in \Omega_f : f(\omega_f, \cdot) \in L^2(\mathcal{D})  \right) = 1. }
			\end{equation}
		\end{enumerate}
	\end{assumption}
	
	We point out that in Assumption~\ref{assum:A-f}~\textit{\ref{assum:matrix-A}}, when $\mathcal{D}$ is a two-dimensional domain ($d=2$), condition \textit{\ref{assum:uniform-ellipticity}} automatically implies \textit{\ref{assum:Cordes-condition}}. However, this does not necessarily hold in higher dimensions ($d \geq 3$).
	
	To transition from the sample space $\Omega$ to a real space, we assume that the input data $\boldsymbol{A}$ and $f$ are represented by random variables. In many applications, the source of randomness can be approximated using a finite number of uncorrelated or independent random variables. A notable example of this is the truncated Karhunen-Loève expansion \cite{Levy}. This motivates our focus on problems involving a finite number of random variables.
	
	\vspace*{1mm}
	\begin{assumption}[parametrization of random inputs]
		\label{assum:random-data}
		The input data $\boldsymbol{A}(\omega_A, \boldsymbol{x})$ and $f(\omega_f, \boldsymbol{x})$ have the form 
		\begin{equation}
			\begin{split}
				\boldsymbol{A}(\omega_A, \boldsymbol{x}) &= 
				\boldsymbol{A}(\boldsymbol{y}_{A}(\omega_A), \boldsymbol{x})
				~~ \mathrm{in}~  \Omega_{A} \times \mathcal{D},
				\\
				f(\omega_f, \boldsymbol{x}) &= 
				f(\boldsymbol{y}_{f}(\omega_f), \boldsymbol{x})
				~~~~ \mathrm{in}~ \Omega_{f} \times \mathcal{D},
			\end{split}
		\end{equation}
		where $\boldsymbol{y}_{A}(\omega_A)= \left( 
		y_{A,1}(\omega_{A}), \cdots, y_{A,N_{A}}(\omega_{A})
		\right) $
		is a vector of real-valued uncorrelated random variables and likewise for 
		$\boldsymbol{y}_{f}(\omega_f)= \left( 
		y_{f,1}(\omega_{f}), \cdots, y_{f,N_{f}}(\omega_{f})
		\right) $
		with $N_{A}, N_{f} \in \mathbb{N}$.
	\end{assumption}
	
	\vspace*{1mm}
	Now we define $\boldsymbol{y} := (\boldsymbol{y}_A, \boldsymbol{y}_f) = (y_1, \cdots, y_N)$, where $N= N_A + N_f$. By considering the random variable $\left\lbrace y_n \right\rbrace _{n=1}^{N}$ which maps the sample space $\Omega$ to the real space $\mathbb{R}^N$, we let $\Gamma_n =y_{n}(\Omega) \subset \mathbb{R} $ denote the image of the random variable $y_n$, and set $\boldsymbol{\Gamma} = \prod_{n= 1}^{N} {\Gamma_n}$. We also set the joint probability density function (PDF) of $\left\lbrace y_n \right\rbrace  _{n= 1}^N$ denoted by $\rho(\boldsymbol{y}):\boldsymbol{\Gamma} \rightarrow \mathbb{R}^{+} \cup \lbrace 0 \rbrace $, with $\rho \in L^\infty(\boldsymbol{\Gamma})$.
	Indeed, in Assumption~\ref{assum:random-data}, the probability space $(\Omega, \mathcal{F}, P)$ is mapped to $(\boldsymbol{\Gamma}, \mathcal{B}(\boldsymbol{\Gamma}), \rho(\boldsymbol{y}) d\boldsymbol{y})$, where $\mathcal{B}(\boldsymbol{\Gamma})$ is the Borel $\sigma$-algebra on $\boldsymbol{\Gamma}$ and $\rho(\boldsymbol{y}) d\boldsymbol{y}$ is a finite measure. 
	For any $\boldsymbol{\phi} \in L^1_{\rho}( \boldsymbol{\Gamma}, \mathbb{R}^k), ~ k \in \mathbb{N} $, we denote its	expected value by $\mathrm{E}_{\rho}[\boldsymbol{\phi}] = \int_{\boldsymbol{\Gamma}} \boldsymbol{\phi}(\boldsymbol{y}) \rho(\boldsymbol{y}) d \boldsymbol{y}$.
	
	Assumption~\ref{assum:random-data} and Doob-Dynkin lemma \cite{Oksendal} guarantee that $u$, the solution of SBVP (\ref{eq:nondivergence}), depends on the same random variables as the input data $\boldsymbol{A}, f$; i.e., $u(\omega, \boldsymbol{x}) = u(y_1(\omega),\dots, y_N(\omega), \boldsymbol{x})$. Then, it is natural to treat $u(\boldsymbol{y}, \boldsymbol{x})$, a function of $N$ random parameters and $d$ spatial variables, as a function of $N+d$ variables.
	
	Through switching to the probability space $(\boldsymbol{\Gamma}, \mathcal{B}(\boldsymbol{\Gamma}), \rho(\boldsymbol{y}) d\boldsymbol{y})$, the problem~(\ref{eq:nondivergence}) turns to finding a stochastic function, $u :\boldsymbol{\Gamma} \times \bar{\mathcal{D}} \rightarrow \mathbb{R} $, such that a.s., the following equation holds
	\begin{equation}
		\label{eq:nondivergence-random-variable}
		\boldsymbol{A}(\boldsymbol{y}, \cdot): \mathrm{D}^2u(\boldsymbol{y}, \cdot) = f( \boldsymbol{y}, \cdot) ~ ~ \mathrm{in}~ \mathcal{D},
		\quad \mathrm{and}~~
		u(\boldsymbol{y}, \cdot) =0 ~~ \mathrm{on}~\partial \mathcal{D}.
	\end{equation}
	
	\vspace*{2mm}
	\subsection{Functional spaces}
	For Hilbert spaces $H_1$ and $H_2$, the Hilbert tensor space $H_1 \otimes H_2$ is defined by the completion of formal sums $u(x, y) = \sum_{i=1}^n{v_i(x)w_i(y)}, ~ \left\lbrace v_i\right\rbrace \subset H_1, ~ \left\lbrace  w_i\right\rbrace  \subset H_2$, with respect to the inner product $\left( u, \hat{u}\right) _{H_1 \otimes H_2}= \sum_{i,j}{(v_i, \hat{v}_j)_{H_1}(w_i, \hat{w}_j)_{H_2}}$ \cite{Babuska-et-al}. Thus, for every $v \in L^2_{\rho}(\boldsymbol{\Gamma}) \otimes \left( H^2(\mathcal{D}) \cap H^1_0(\mathcal{D})\right) $, we have 
	\begin{equation*}
		\lVert v  \rVert_{L^2_{\rho}(\boldsymbol{\Gamma}) \otimes H^2(\mathcal{D})}^2: = \int_{\boldsymbol{\Gamma}} {\lVert v( \boldsymbol{y}, \cdot)  \rVert_{H^2(\mathcal{D})}^2 \rho(\boldsymbol{y})  d\boldsymbol{y}} = \mathrm{E}_{\rho} [\lVert  v(\boldsymbol{y}, \cdot) \rVert_{H^2(\mathcal{D})}^2]. 
	\end{equation*}
	Indeed, $v( \boldsymbol{y}, \cdot) \in  H^2(\mathcal{D}) \cap H^1_0(\mathcal{D}) $ for a.e. $\boldsymbol{y} \in \boldsymbol{\Gamma}$ and $v(\cdot, \boldsymbol{x}) \in L^2_{\rho}(\boldsymbol{\Gamma})$ for a.e. $\boldsymbol{x} \in \mathcal{D}$ and moreover we have the isomorphism 
	$L^2_{\rho}(\boldsymbol{\Gamma}) \otimes \left( H^2(\mathcal{D}) \cap H^1_0(\mathcal{D})\right)  \cong  L^2_{\rho}(\boldsymbol{\Gamma}; H^2(\mathcal{D}) \cap H^1_0(\mathcal{D}))$.
	
	We introduce some notations that will be used throughout this paper. Additional notations will be defined as new spaces are introduced. The following are the notations for some function spaces
	\begin{equation*}
		\begin{aligned}
			&H ~:= H^1_0(\mathcal{D}), \quad
			&&LH ~:= L^2_{\rho}(\boldsymbol{\Gamma}) \otimes H^1_0(\mathcal{D}), \\
			&\boldsymbol{H} ~:= H^1(\mathcal{D};\mathbb{R}^d), \quad
			&&L\boldsymbol{H} ~:= L^2_{\rho}(\boldsymbol{\Gamma}) \otimes H^1(\mathcal{D};\mathbb{R}^d), \\
			&H^2 := H^2(\mathcal{D}) \cap H^1_0(\mathcal{D}), \quad
			&&LH^2 := L^2_{\rho}(\boldsymbol{\Gamma}) \otimes (H^2(\mathcal{D}) \cap H^1_0(\mathcal{D})).
		\end{aligned}
	\end{equation*}
	All derivative-based operators considered in the following, including the Laplacian $\Delta$, gradient $\nabla$, derivative $\mathrm{D}$,  divergence $\nabla \cdot$, and curl $\nabla \times$, are taken with respect to the spatial variable $\boldsymbol{x}$.

	\vspace*{1mm}
	\begin{lemma}[well-posedness of the stochastic Laplace equation {\cite[Section 2.2]{Babuska-et-al}}]
		\label{rem:well-posedness-laplac-eq}
		The problem of finding a strong solution $u \in LH^2$ to the equation
		\begin{equation}
			\label{eq:laplac}
			\Delta u(\boldsymbol{y}, \cdot) = f (\boldsymbol{y}, \cdot) 
			~~ \mathrm{in}~\mathcal{D},
			~~~ \mathrm{and}~~
			u(\boldsymbol{y}, \cdot) = 0 ~ ~\mathrm{on}~ \partial \mathcal{D},
		\end{equation}
		for any given $f \in L^2_{\rho}(\boldsymbol{\Gamma}, L^2(\mathcal{D}))$ is well-posed. It means that there exists a unique solution $u(\boldsymbol{y}, \cdot) \in H^2$ that satisfies (\ref{eq:laplac}) for a.e. $\boldsymbol{y} \in \boldsymbol{\Gamma}$; moreover, there exists $C_{\ref{ineq:stability-laplac}}>0$ such that for every $ v \in LH^2$ the following holds
		\begin{equation}
			\label{ineq:stability-laplac}
			\rVert v \lVert_{L^2_{\rho}(\boldsymbol{\Gamma}) \otimes H^2(\mathcal{D})}^2
			=
			\mathrm{E}_{\rho}\left[ \lVert v\rVert_{H^2(\mathcal{D})}^2 \right] 
			\leq 
			C_{\ref{ineq:stability-laplac}}
			\mathrm{E}_{\rho} \left[ \lVert \Delta v \rVert _{L^2(\mathcal{D})}^2 \right]  
			= C_{\ref{ineq:stability-laplac}} 
			\rVert  \Delta v \rVert _{L^2_{\rho}(\boldsymbol{\Gamma}) \otimes L^2(\mathcal{D})}^2
			.
		\end{equation}
	\end{lemma}
	
	\vspace*{1mm}
	\begin{definition}[strong solution] 
		The function $u \in LH^2 $ is a strong solution to problem (\ref{eq:nondivergence-random-variable}), if 
		for a.e. $\boldsymbol{y} \in \boldsymbol{\Gamma}$, $u(\boldsymbol{y}, \cdot)$ satisfies (\ref{eq:nondivergence-random-variable}) a.e.;  equivalently when the following holds 
		\begin{multline}
			\label{eq:weak-form}
			\int_{\boldsymbol{\Gamma}} \int_{\mathcal{D}} {\left[ \boldsymbol{A}( \boldsymbol{y}, \boldsymbol{x}): \mathrm{D}^2u( \boldsymbol{y}, \boldsymbol{x})\right]  \Delta v(\boldsymbol{y}, \boldsymbol{x}) \rho(\boldsymbol{y}) d\boldsymbol{x}  d \boldsymbol{y} }  
			\\
			=
			\int_{\boldsymbol{\Gamma}} \int_{\mathcal{D}} { f(\boldsymbol{y}, \boldsymbol{x}) \Delta v(\boldsymbol{y}, \boldsymbol{x}) \rho(\boldsymbol{y}) d\boldsymbol{x}  d \boldsymbol{y} } 
			\quad \forall v \in LH^2.
		\end{multline}
	\end{definition}
	
	\vspace*{2mm}
	\subsection{Existence and uniqueness of the strong solution}
	To establish the well-posedness of the strong solution, it is convenient, as suggested by the approach in \cite{Smears-Suli}, to first prove the well-posedness of the strong solution to an appropriately scaled equation, such as $\gamma\boldsymbol{A}:\mathrm{D}^2 u = \gamma f.$ 
	In this context, we define the stochastic scaling function $\gamma:= \frac{\mathrm{tr} \boldsymbol{A}}{\lvert \boldsymbol{A} \rvert^2}$. From Assumption~\ref{assum:A-f}~\textit{\ref{assum:matrix-A}\ref{assum:uniform-ellipticity}}, we deduce that there exists a constant $\Lambda \in (0, 1]$ such that for a.e. $\boldsymbol{y} \in \boldsymbol{\Gamma}$, the following holds
	\begin{equation}
		\label{set:boundedness:gamma} 
		\Lambda \leq
		\lVert \gamma(\boldsymbol{y}, \cdot) \rVert_{L^{\infty}(\mathcal{D})}
		\leq
		\Lambda^{-1}.
	\end{equation}
	Define the bilinear form $a_{\gamma}: LH^2 \times LH^2  \rightarrow \mathbb{R}$ by
	\begin{multline}
		\label{def:bilinear-form}
		a_{\gamma}(v, w):= \int_{\boldsymbol{\Gamma}} \int_{\mathcal{D}} {\left[ \gamma 	\boldsymbol{A}(\boldsymbol{y}, \boldsymbol{x}): \mathrm{D}^2v(\boldsymbol{y}, \boldsymbol{x})\right]  \Delta w(\boldsymbol{y}, \boldsymbol{x}) \rho(\boldsymbol{y})d\boldsymbol{x}  d \boldsymbol{y} } 
		\\
		=
		\mathrm{E}_{\rho}\left[  
		\int_{\mathcal{D}}{\left[ \gamma \boldsymbol{A}: \mathrm{D}^2 v\right] \Delta w  d\boldsymbol{x}}
		\right] .
	\end{multline} 
	Also define the linear functional $F_{\gamma}: LH^2 \rightarrow \mathbb{R}$ by 
	\begin{equation}
		\label{def:right-hand-side-functional}
		F_{\gamma}(v):= \int_{\boldsymbol{\Gamma}} \int_{\mathcal{D}} { \gamma f(\boldsymbol{y}, 		\boldsymbol{x}) \Delta v(\boldsymbol{y}, \boldsymbol{x}) \rho(\boldsymbol{y})  d\boldsymbol{x}  d \boldsymbol{y} 
		}
		=  
		\mathrm{E}_{\rho} \left[ \int_{\mathcal{D}}{ \gamma f \Delta vd \boldsymbol{x}} 
		\right].
	\end{equation}

		\vspace*{1mm}
	\begin{lemma}[well-posedness]
		\label{lem:wellposedness-L-M}
		The problem of finding a strong solution $u \in LH^2$ such that
		\begin{equation}
			\label{eq:scaled-bilinear-formed-problem}
			a_{\gamma}(u, v)= F_{\gamma}(v) \quad \forall v \in LH^2,
		\end{equation}
		is well-posed.
	\end{lemma}
	
	\vspace*{1mm}
	\begin{proof}
		%		A direct application of the Lax-Milgram theorem ensures the existence and uniqueness of the solution.
		The argument follows the same reasoning as in Theorem 3 of \cite{Smears-Suli}, with the only modification being the replacement of $LH^2$ in place of $H^2$.
	\end{proof}
	
	\vspace*{1mm}
	\begin{remark}[equivalence of scaled and main equations]
		The strict positivity of the scaling factor $\gamma$, as concluded from Assumption~\ref{assum:A-f}~\textit{\ref{assum:matrix-A}\ref{assum:uniform-ellipticity}}, together with the uniform boundedness in (\ref{set:boundedness:gamma}), implies that the solutions of (\ref{eq:weak-form}) and  (\ref{eq:scaled-bilinear-formed-problem}) are identical. Consequently, Lemma~\ref{lem:wellposedness-L-M} ensures that the problem of finding a strong solution to (\ref{eq:weak-form}) is also well-posed.
	\end{remark}	
	
	\vspace*{2mm}
	\section{Mixed formulation}
	\label{sec:Mixed-formulation}
	
	Since numerical approximation requires seeking the solution of (\ref{eq:weak-form}) in a finite dimensional subspace of the tensor product space $LH^2$, dealing with such a large regular function space ($H^2(\mathcal{D})$) in the physical domain $\mathcal{D}$ can lead to complicated and unpleasant computations. Therefore, as suggested in \cite{Gallistl}, we propose considering an alternative, equivalent problem in which the solution resides in a weaker space, as discussed in \cite{Lakkis-Mousavi}. 
	In this regard, we define the space
	\begin{equation*}
		\boldsymbol{H}_{\boldsymbol{t}}:= \left\lbrace  \boldsymbol{\psi} \in H^1(\mathcal{D}; \mathbb{R}^d) {\big \vert}  ~  \boldsymbol{\psi}_{\boldsymbol{t}} := \boldsymbol{\psi} - (\boldsymbol{\psi}\cdot \boldsymbol{n})\boldsymbol{n} = \boldsymbol{0} \right\rbrace,
	\end{equation*}
	where $\boldsymbol{n}(\boldsymbol{x})$ denotes the outward unit normal vector to $\mathcal{D}$ for $\mathcal{S}$-almost every $\boldsymbol{x} \in \partial \mathcal{D}$. Hence, $\boldsymbol{\psi}_{\boldsymbol{t}}$ represents the tangential component of $\boldsymbol{\psi}$ on the boundary 
	$\partial \mathcal{D}$. We also introduce the notation 
	\begin{equation*}
		L\boldsymbol{H}_{\boldsymbol{t}}:= L^2_{\rho}(\boldsymbol{\Gamma}) \otimes \boldsymbol{H}_{\boldsymbol{t}}.
	\end{equation*}
	Using this, we define the following stochastic quadratic functional on $LH \times L\boldsymbol{H}_{\boldsymbol{t}}$
	\begin{equation}
		\label{fun:cost-functional}
		E(\varphi, \boldsymbol{\psi}):= \lVert \nabla \varphi - \boldsymbol{\psi} \rVert^2_{ L^2_\rho(\boldsymbol{\Gamma}) \otimes L^2(\mathcal{D})}
		+
		\lVert \nabla \times \boldsymbol{\psi} \rVert^2_{ L^2_\rho(\boldsymbol{\Gamma}) \otimes L^2(\mathcal{D})} 
		+
		\lVert \boldsymbol{A} : \mathrm{D} \boldsymbol{\psi} - f\rVert_{ L^2_\rho(\boldsymbol{\Gamma}) \otimes L^2(\mathcal{D})}^2.
	\end{equation}
	Here, $\nabla \times $ denotes the curl operator
	\begin{equation*}
		\nabla \times : H^1(\mathcal{D}; \mathbb{R}^d) \rightarrow L^2(\mathcal{D})^{\hat{d}} \quad \text{for  } \hat{d}:= \binom{d}{2}= \begin{cases}
			1 & \text{ if } d=2,
			\\
			3 & \text{ if } d=3.
		\end{cases}
	\end{equation*} 
	The curl operator for two- and three-component vector fields, respectively, is defined as
	\begin{equation}
		\label{eq:curl-operator}
		\nabla \times \\
		\begin{bmatrix}
			\psi_1
			\\
			\psi_2
		\end{bmatrix}
		=
		\partial_{x_1} \psi_2 - \partial_{x_2} \psi_1,
		\qquad
		\nabla \times \\
		\begin{bmatrix}
			\psi_1
			\\
			\psi_2
			\\
			\psi_3
		\end{bmatrix}
		= 
		\begin{bmatrix}
			\partial_{x_2} \psi_3 - \partial_{x_3} \psi_2
			\\
			\partial_{x_3} \psi_1 - \partial_{x_1} \psi_3
			\\
			\partial_{x_1} \psi_2 - \partial_{x_2} \psi_1
		\end{bmatrix}.
	\end{equation}
	We then deal with the convex minimization problem of finding a unique pair of the form
	\begin{equation}
		\label{eq:minimization}
		(u, \boldsymbol{g})= 
		\underset
		{\substack{
				(\varphi, \boldsymbol{\psi})\in 
				LH \times L\boldsymbol{H}_{\boldsymbol{t}}
		}}
		{\mathrm{argmin}}
		E(\varphi ,\boldsymbol{\psi}).
	\end{equation}
	
	\vspace*{1mm}
	\begin{remark}[equivalence of (\ref{eq:nondivergence-random-variable}) and (\ref{eq:minimization})]
		\label{rem:equivalent-problems}
		If $u \in LH^2$ is a strong solution to (\ref{eq:nondivergence-random-variable}), then $(u, \nabla u)$ minimizes the non-negative convex functional $E$, causing it to attain a zero value.
		Conversely, if $E$ attains its minimum value at $(u, \boldsymbol{g})$, then $u$ is also a strong solution to (\ref{eq:nondivergence-random-variable}) and	$ \boldsymbol{g} = \nabla u$, in $L^2_{\rho}(\boldsymbol{\Gamma}) \otimes L^2(\mathcal{D})$. 
		Therefore the problem of finding a strong solution to (\ref{eq:nondivergence-random-variable}) and problem (\ref{eq:minimization}) are equivalent. Throughout the rest of the paper, the symbol $\boldsymbol{g}$ will be used interchangeably with $\nabla u$.
	\end{remark}
	
	\vspace*{1mm}
	The $L^2_\rho(\boldsymbol{\Gamma}) \otimes L^2(\mathcal{D})$, $L^2_\rho(\boldsymbol{\Gamma}) \otimes L^2(\partial \mathcal{D})$, $L^2(\mathcal{D})$ and $L^2(\partial \mathcal{D})$, inner products of two scalar, vector, or tensor-valued functions $\varphi$ and $\psi$ are indicated with
	\begin{equation*}
		\begin{aligned}
			&\left\langle \varphi, \xi \right\rangle_{\rho, \mathcal{D}} := \int_{\boldsymbol{\Gamma}} \int_{\mathcal{D}} \hspace*{-0.5mm}
			\varphi(\boldsymbol{y}, \boldsymbol{x}) \star \psi (\boldsymbol{y}, \boldsymbol{x}) \rho(\boldsymbol{y})  d\boldsymbol{x}
			 d\boldsymbol{y},
			&&\left\langle \varphi, \xi \right\rangle_{\mathcal{D}} :=  \int_{\mathcal{D}} 
			\hspace*{-0.5mm}
			\varphi (\boldsymbol{x}) \star \psi (\boldsymbol{x})  d\boldsymbol{x},
			\\
			&\left\langle \varphi, \xi \right\rangle_{\rho,\partial \mathcal{D}} := \int_{\boldsymbol{\Gamma}} \int_{\partial \mathcal{D}} \hspace*{-2mm} \varphi (\boldsymbol{y}, \boldsymbol{x}) \star \psi (\boldsymbol{y}, \boldsymbol{x}) \rho(\boldsymbol{y}) d\mathcal{S} (\boldsymbol{x}) d\boldsymbol{y},
			%\hspace*{-2mm}
			&&
			\left\langle \varphi, \xi \right\rangle_{\partial \mathcal{D}} := \int_{\partial \mathcal{D}} \hspace*{-2mm} \varphi (\boldsymbol{x}) \star \psi (\boldsymbol{x}) d\mathcal{S} (\boldsymbol{x}),
		\end{aligned}
	\end{equation*}
	where $\star$ stands for one of the arithmetic, Euclidean-scalar, or Frobenius inner product in $\mathbb{R}$, $\mathbb{R}^d$, or $\mathbb{R}^{d \times d}$ respectively. 
	For any $(\varphi, \boldsymbol{\psi}) \in LH \times L\boldsymbol{H}$ we define 
	\begin{equation}
		\label{eq:H1-norm}
		\lVert  (\varphi, \boldsymbol{\psi}) \rVert^2_{L^2_\rho(\boldsymbol{\Gamma}) \otimes  H^1(\mathcal{D})} := 
		\lVert  \varphi \rVert^2_{L^2_\rho(\boldsymbol{\Gamma}) \otimes  H^1_0(\mathcal{D})} 
		+
		\lVert  \boldsymbol{\psi} \rVert^2_{L^2_\rho(\boldsymbol{\Gamma}) \otimes  H^1(\mathcal{D})}.
	\end{equation} 
	The Euler–Lagrange equation of the minimization problem (\ref{eq:minimization}) consists in finding $(u, \boldsymbol{g}) \in LH \times L\boldsymbol{H}_{\boldsymbol{t}} $ such that
	\begin{multline}
		\label{eq:Euler-Lagrange-stochastic}
		\left\langle \nabla u - \boldsymbol{g}, \nabla\varphi - \boldsymbol{\psi} \right\rangle_{\rho, \mathcal{D}}  
		+
		\left\langle  \nabla \times \boldsymbol{g}, \nabla \times \boldsymbol{\psi} \right\rangle_{\rho, \mathcal{D}} 
		+
		\left\langle  \boldsymbol{A} : \mathrm{D} \boldsymbol{g}, \boldsymbol{A} : \mathrm{D} \boldsymbol{\psi} \right\rangle_{\rho, \mathcal{D}} 
		\\
		= 
		\left\langle  f, \boldsymbol{A} : \mathrm{D} \boldsymbol{\psi} \right\rangle_{\rho, \mathcal{D}} 
		~~
		\forall (\varphi, \boldsymbol{\psi}) \in LH \times L\boldsymbol{H}_{\boldsymbol{t}}  .
	\end{multline}
	In accordance with (\ref{eq:Euler-Lagrange-stochastic}), we define the stochastic symmetric bilinear form $a: \left( LH \times L\boldsymbol{H}_{\boldsymbol{t}}  \right)^2  \rightarrow \mathbb{R}$ by
	%	\begin{equation*}
		%		a: \left( LH \times L\boldsymbol{H}_t  \right)^2  \rightarrow \mathbb{R}
		%	\end{equation*}
	%	by
	\begin{equation}
		\label{def:stochastic-bilinear-a}
		a (\varphi, \boldsymbol{\psi}; \varphi^{\prime}, \boldsymbol{\psi}^{\prime}):=
		\left\langle \nabla\varphi - \boldsymbol{\psi} , \nabla\varphi^{\prime} - 	\boldsymbol{\psi}^{\prime}  \right\rangle_{\rho, \mathcal{D}}  
		+
		\left\langle  \nabla \times \boldsymbol{\psi}, \nabla \times 	\boldsymbol{\psi}^{\prime} \right\rangle_{\rho, \mathcal{D}} 
		+
		\left\langle  \boldsymbol{A} : \mathrm{D} \boldsymbol{\psi}, \boldsymbol{A} : \mathrm{D} 	\boldsymbol{\psi}^{\prime} \right\rangle_{\rho, \mathcal{D}}. 
	\end{equation}
	
	\vspace*{1mm}
	\begin{lemma}[coercivity and continuity of $a(\cdot; \cdot)$]
		\label{lem:well-posedness-a-y}
		The stochastic bilinear form $a$ is coercive and continuous; i.e., there exist $C_{\ref{ineq:coercivity-a}}, C_{\ref{ineq:continuity-a}}> 0$ such that for 
		any $(\varphi, \boldsymbol{\psi}), (\varphi^{\prime}, \boldsymbol{\psi}^{\prime}) \in LH \times L\boldsymbol{H}_{\boldsymbol{t}}$ the following hold
		\begin{multline}
			\label{ineq:coercivity-a}
			a (\varphi, \boldsymbol{\psi}; \varphi, \boldsymbol{\psi}) = 
			\lVert \nabla \varphi - \boldsymbol{\psi} \rVert^2_{L^2_\rho(\boldsymbol{\Gamma}) \otimes L^2(\mathcal{D})}
			+
			\lVert \nabla \times \boldsymbol{\psi} \rVert^2_{L^2_\rho(\boldsymbol{\Gamma}) \otimes L^2(\mathcal{D})} 
			+
			\lVert \boldsymbol{A} : \mathrm{D}  \boldsymbol{\psi} \rVert_{L^2_\rho(\boldsymbol{\Gamma}) \otimes L^2(\mathcal{D})}^2
			\\
			\geq
			C_{\ref{ineq:coercivity-a}}
			\left\| (\varphi, \boldsymbol{\psi})  \right\| ^2_{ L^2_\rho(\boldsymbol{\Gamma}) \otimes H^1(\mathcal{D}) }
			,
		\end{multline}
		\\
		\vspace*{-1.2cm}
		\begin{equation}
			\label{ineq:continuity-a}
			%\hspace*{-3.8cm}
			\left|  a (\varphi, \boldsymbol{\psi}; \varphi^\prime, \boldsymbol{\psi})^\prime \right| 
			\leq
			C_{\ref{ineq:continuity-a}}
			\left\|  (\varphi, \boldsymbol{\psi}) \right\|_{
				L^2_\rho(\boldsymbol{\Gamma}) \otimes H^1(\mathcal{D}) } 
			\left\|  (\varphi^{\prime}, \boldsymbol{\psi}^\prime) \right\|_{   L^2_\rho(\boldsymbol{\Gamma}) \otimes H^1(\mathcal{D})
			}.
		\end{equation}
	\end{lemma}
	
	\vspace*{1mm}
	\begin{proof}
		The argument to demonstrate coercivity follows the same reasoning as in Theorem 3.7 of \cite{Lakkis-Mousavi}, while the argument for continuity is based on Section 3.9 of the same reference. The only modification involves replacing $LH \times L\boldsymbol{H}_{\boldsymbol{t}}$ in place of $H \times \boldsymbol{H}_{\boldsymbol{t}}$.
	\end{proof}
	
	\vspace*{1mm}
	From another perspective, we can consider the solution $(u, \boldsymbol{g})$ as functions $u:\boldsymbol{\Gamma}\rightarrow H$ and $\boldsymbol{g}: \boldsymbol{\Gamma}\rightarrow \boldsymbol{H}_{\boldsymbol{t}}$ respectively. We use the notation $u(\boldsymbol{y}), \boldsymbol{g}(\boldsymbol{y}),  \boldsymbol{A}(\boldsymbol{y}), f(\boldsymbol{y})$ whenever we want to emphasize on the dependence on parameter $\boldsymbol{y}$. Therefore the problem (\ref{eq:nondivergence-random-variable}) is equivalent to finding $\left( u(\boldsymbol{y}), \boldsymbol{g}(\boldsymbol{y}) \right) \in H \times \boldsymbol{H}_{\boldsymbol{t}}$ such that the following equation holds
	\begin{multline}
		\label{eq:Euler-Lagrange-parametric}
		\left\langle \nabla u (\boldsymbol{y}) - \boldsymbol{g} (\boldsymbol{y}), \nabla\varphi - \boldsymbol{\psi} \right\rangle_{\mathcal{D}}  
		+
		\left\langle  \nabla \times \boldsymbol{g}(\boldsymbol{y}), \nabla \times \boldsymbol{\psi} \right\rangle_{\mathcal{D}} 
		\\
		+
		\left\langle  \boldsymbol{A}(\boldsymbol{y}):\mathrm{D} \boldsymbol{g}(\boldsymbol{y}), \boldsymbol{A}(\boldsymbol{y}):\mathrm{D} \boldsymbol{\psi} \right\rangle_{\mathcal{D}} 
		%\\
		= 
		\left\langle  f(\boldsymbol{y}), \boldsymbol{A}(\boldsymbol{y}):\mathrm{D} \boldsymbol{\psi} \right\rangle_{\mathcal{D}} 
		~~
		\forall (\varphi, \boldsymbol{\psi}) \in 
		H \times \boldsymbol{H}_{\boldsymbol{t}} .
	\end{multline}
	In this case, we define the parametric bilinear form $	a_{\boldsymbol{y}}: \left( H \times \boldsymbol{H}_{\boldsymbol{t}} \right)^2  \rightarrow \mathbb{R}$ as
	%		\begin{equation*}
		%			a_{\boldsymbol{y}}: \left( H \times \boldsymbol{H}_t \right)^2  \rightarrow \mathbb{R}
		%		\end{equation*}
	%		as
	\begin{multline}
		\label{def:parametric-bilinear-a}
		a_{\boldsymbol{y}} (\varphi, \boldsymbol{\psi}; \varphi^{\prime}, \boldsymbol{\psi}^{\prime}):=
		\left\langle \nabla\varphi - \boldsymbol{\psi} , \nabla\varphi^{\prime} - \boldsymbol{\psi}^{\prime}  \right\rangle_{\mathcal{D}}  
		+
		\left\langle  \nabla \times \boldsymbol{\psi}, \nabla \times \boldsymbol{\psi}^{\prime} \right\rangle_{\mathcal{D}} 
		\\
		+
		\left\langle  \boldsymbol{A}(\boldsymbol{y}):\mathrm{D}\boldsymbol{\psi}, \boldsymbol{A}(\boldsymbol{y}):\mathrm{D}\boldsymbol{\psi}^{\prime} \right\rangle_{\mathcal{D}}. 
	\end{multline}
	
	\vspace*{2mm}
	\section{Discretization}
	\label{sec:discretization}
	
	In this section, we discuss the discretization of the solution in the physical and stochastic domains separately. We begin with the continuous finite element method for spatial discretization and then proceed to the collocation method for stochastic discretization.
	
	\vspace*{2mm}
	\subsection{Continuous finite element method}
	\label{susec:continuous-FEM}
	Consider $\mathfrak{T}$ as a collection of shape regular conforming simplicial partitions, also known as triangulations, of $\mathcal{D}$ into simplices. For a given $\mathcal{T} \in \mathfrak{T}$, and for each $K \in \mathcal{T}$ , let $h_K:= \text{diam} K$ and define $h:= \max_{K\in \mathcal{T}} h_K$ . Let $\mathcal{F}^b$ represent the set of boundary faces, $\cup_{F \in \mathcal{F}^b} = \partial \mathcal{D}$, and for each $F \in \mathcal{F}^b$, let $h_{K,F}$ denote the diameter of the element associated with the boundary face $F$. 
	When the faces of elements are chosen to be flat, a curved boundary $\partial \mathcal{D}$ prevents $\cup_{K \in \mathcal{T}} K$ from coinciding exactly with $\mathcal{D}$. In such cases, sections of $\partial \mathcal{D}$ can be approximated using line segments or elements with curved sides for the corresponding boundary elements. For simplicity, we assume that $\mathcal{D}$ is a polytopal domain for the remainder of this paper. However, the results remain valid for a general convex domain.
	
	Enforcing the vanishing tangential trace constraint in finite element spaces is not computationally straightforward. Additionally, if we attempt to enforce this constraint by adding a squared term to the cost functional $E$, establishing the coercivity of the corresponding bilinear form becomes challenging. Therefore, in this work, by utilizing an inverse trace inequality in finite element spaces, we introduce a discrete cost functional and handle this issue within a discrete functional framework. The Galerkin finite element spaces are defined as
	\begin{equation}
		\label{def:Galerkin-spaces}
		\begin{aligned}
		\mathbb{U} &:= 
		\left\lbrace v \in H^1_0(\mathcal{D}),~ v|_{K} \in  
		\mathcal{P}_k(K), ~ \forall K \in \mathcal{T} \right\rbrace, 
		\\
		\mathbb{G} &:= 
		 \left\lbrace \boldsymbol{v} \in H^1(\mathcal{D}; \mathbb{R}^d) , ~ 
		 \boldsymbol{v}|_{K} \in  
		 \mathcal{P}_k(K; \mathbb{R}^d), ~ \forall K \in \mathcal{T}, \right\rbrace.
		\end{aligned}
	\end{equation}
%
%	\begin{equation}
%		\label{def:Galerkin-spaces}
%		\mathbb{U}:= \mathbb{P}_k(\mathcal{T}) \cap H^1_0(\mathcal{D}),
%		\quad
%		\mathbb{G}:= \mathbb{P}_k(\mathcal{T}; \mathbb{R}^d) \cap H^1(\mathcal{D}; \mathbb{R}^d),
%	\end{equation}
%	Here, $\mathcal{P}_k(K)$ and $\mathcal{P}_k(K; \mathbb{R}^d)$ denote the scalar-valued and vector-valued spaces of polynomials of total degree $k$ or less on the element $K$, respectively.
	Here, $\mathcal{P}_k(K)$ and $\mathcal{P}_k(K; \mathbb{R}^d)$ denote the scalar-valued and vector-valued polynomial spaces of degree at most $k \in \mathbb{N}$ on the element $K$, respectively.
	The notation $N_{h,k}$ denotes the dimension of the finite element space $\mathbb{U} \times \mathbb{G}$.
%	
%	 represent the scalar valued polynomial of degree $k$ on $K$ and the vector valued polynomial of degree $k$ on $K$ respectively and 
%	and $N_{h,k}$ denotes the dimension of $\mathbb{U} \times \mathbb{G}$ finite element space.

	\vspace*{1mm}
	\begin{definition} [mesh-dependent $H^1$-norm]
		\label{def:discrete-norm}
		We define the mesh-dependent $H^1$-type norm, first on $L\boldsymbol{H}$
		as 
		%\begin{multline}
		\begin{equation}
			\label{eq:discrete-norm-G}
			\begin{aligned}
			\lVert \boldsymbol{\psi} \rVert^2_{L^2_{\rho}(\boldsymbol{\Gamma}) \otimes H^1_h(\mathcal{D})}:= 
			&\lVert \mathrm{D} \boldsymbol{\psi}  \rVert_{L^2_{\rho}(\boldsymbol{\Gamma}) \otimes L^2(\mathcal{D})}^2 
			+
			\lVert \boldsymbol{\psi}_{\boldsymbol{t}} \rVert^2_{L^2_{\rho}(\boldsymbol{\Gamma}) \otimes L^2(\partial \mathcal{D})}
			\\
			&+
			\sum_{F \in \mathcal{F}^b}
			h_{K,F}^{-1}\lVert \boldsymbol{\psi}_{\boldsymbol{t}} \rVert^2_{L^2_{\rho}(\boldsymbol{\Gamma}) \otimes L^2(F)},
		\end{aligned}
		\end{equation}
		%\end{multline}
		and then on $ LH \times L \boldsymbol{H}$ as 
		\begin{equation}
			\label{eq:discrete-norm-U-times-G}
			\lVert (\varphi, \boldsymbol{\psi})  \rVert_{L^2_{\rho}(\boldsymbol{\Gamma}) \otimes H^1_h( \mathcal{D})}^2
			:=
			\lVert \varphi \rVert_{L^2_{\rho}(\boldsymbol{\Gamma}) \otimes H^1_0( \mathcal{D})}^2
			+
			\lVert \boldsymbol{\psi} \rVert^2_{L^2_{\rho}(\boldsymbol{\Gamma}) \otimes H^1_h( \mathcal{D})}.
		\end{equation} 
	\end{definition}
	
	\vspace*{1mm}
	According to (\ref{def:Galerkin-spaces}), we introduce the notations
	\begin{equation*}
		L\mathbb{U}:= L^2_{\rho}(\boldsymbol{\Gamma}) \otimes \mathbb{U}, 
		\quad
		L\mathbb{G}:= L^2_{\rho}(\boldsymbol{\Gamma}) \otimes \mathbb{G}.
	\end{equation*}
	We then define the following mesh-dependent quadratic functional on $L\mathbb{U} \times L\mathbb{G}$
	\begin{multline}
		\label{fun:cost-functional-discerete-local}
		E_h(\varphi, \boldsymbol{\psi}):= \lVert \nabla \varphi - \boldsymbol{\psi} \rVert^2_{L^2_{\rho}(\boldsymbol{\Gamma}) \otimes L^2(\mathcal{D})}
		+
		\lVert \nabla \times \boldsymbol{\psi} \rVert^2_{ L^2_{\rho}(\boldsymbol{\Gamma}) \otimes L^2(\mathcal{D})} 
		\\
		+
		\lVert \boldsymbol{\psi}_{\boldsymbol{t}} \rVert^2_{ L^2_{\rho}(\boldsymbol{\Gamma}) \otimes L^2(\partial \mathcal{D})}
		%\\
		+
		\sum_{F \in \mathcal{F}^b}
		h_{K,F}^{-1}\lVert \boldsymbol{\psi}_{\boldsymbol{t}} \rVert^2_{L^2_{\rho}(\boldsymbol{\Gamma}) \otimes L^2(F)}
		+
		\lVert \boldsymbol{A}:\mathrm{D} \boldsymbol{\psi} - f\rVert_{ L^2_{\rho}(\boldsymbol{\Gamma}) \otimes L^2(\mathcal{D})}^2.
	\end{multline}
	Next, we address the convex minimization problem of finding a unique pair of the form
	\begin{equation}
		\label{eq:minimization-discrete}
		(u_{\mathbb{U}}, \boldsymbol{g}_{\mathbb{G}})= 
		\underset
		{\substack{
				(\varphi, \boldsymbol{\psi})\in 
				L \mathbb{U}  \times L \mathbb{G}
		}}
		{\mathrm{argmin}}
		E_h(\varphi ,\boldsymbol{\psi}).
	\end{equation}
	The Euler–Lagrange equation of the minimization problem (\ref{eq:minimization-discrete}) consists in finding $(u_\mathbb{U}, \boldsymbol{g}_\mathbb{G}) \in  L\mathbb{U} \times L\mathbb{G} $ such that
	\begin{multline}
		\label{eq:Euler-Lagrange-stochastic-discrete-local}
		\left\langle \nabla u_{\mathbb{U}} - \boldsymbol{g}_{\mathbb{G}}, \nabla\varphi - \boldsymbol{\psi} \right\rangle_{\rho,\mathcal{D}}  
		+
		\left\langle  \nabla \times \boldsymbol{g}_{\mathbb{G}}, \nabla \times \boldsymbol{\psi} \right\rangle_{\rho, \mathcal{D}} 
		+
		\left\langle  {\boldsymbol{g}_{\mathbb{G}}}_{\boldsymbol{t}} , \boldsymbol{\psi}_{\boldsymbol{t}} \right\rangle _{\rho,\partial \mathcal{D}}
		\\
		+
		\sum_{F \in \mathcal{F}^b} h_{K,F}^{-1}\left\langle {\boldsymbol{g}_{\mathbb{G}}}_{\boldsymbol{t}} , \boldsymbol{\psi}_{\boldsymbol{t}} \right\rangle _{\rho, F}
		%\\
		+
		\left\langle  \boldsymbol{A}:\mathrm{D} \boldsymbol{g}_{\mathbb{G}}, \boldsymbol{A}:\mathrm{D} \boldsymbol{\psi} \right\rangle_{\rho, \mathcal{D}} 
		\\
		= 
		\left\langle  f, \boldsymbol{A}:\mathrm{D} \boldsymbol{\psi} \right\rangle_{\rho, \mathcal{D}} 
		~~
		\forall (\varphi, \boldsymbol{\psi}) \in 
		L \mathbb{U} \times L \mathbb{G}.
	\end{multline}
	Corresponding to (\ref{eq:Euler-Lagrange-parametric}), we may consider the parametric, $\boldsymbol{y}$-dependent, counterpart of (\ref{eq:Euler-Lagrange-stochastic-discrete-local}) as finding  $(u_\mathbb{U}(\boldsymbol{y}), \boldsymbol{g}_\mathbb{G}(\boldsymbol{y})) \in \mathbb{U} \times \mathbb{G} $ such that
	\begin{multline}
		\label{eq:Euler-Lagrange-stochastic-discrete-local-parametric}
		\left\langle \nabla u_{\mathbb{U}} (\boldsymbol{y}) - \boldsymbol{g}_{\mathbb{G}} (\boldsymbol{y}), \nabla\varphi - \boldsymbol{\psi} \right\rangle_{\mathcal{D}}  
		+
		\left\langle  \nabla \times \boldsymbol{g}_{\mathbb{G}}(\boldsymbol{y}), \nabla \times \boldsymbol{\psi} \right\rangle_{ \mathcal{D}} 
		%+
		%\left\langle  \boldsymbol{t}_D \boldsymbol{g}_{\mathbb{G}} , \boldsymbol{t}_D \boldsymbol{\psi} \right\rangle _{\rho,\partial D}
		\\
		+
		\sum_{F \in \mathcal{F}^b} \left( 1 + h_{K,F}^{-1} \right) \left\langle  \left( \boldsymbol{g}_{\mathbb{G}}(\boldsymbol{y})\right) _{\boldsymbol{t}} , \boldsymbol{\psi}_{\boldsymbol{t}} \right\rangle _{ F}
		+
		\left\langle  \boldsymbol{A}(\boldsymbol{y}):\mathrm{D} \boldsymbol{g}_{\mathbb{G}}(\boldsymbol{y}), \boldsymbol{A}(\boldsymbol{y}):\mathrm{D} \boldsymbol{\psi} \right\rangle_{\mathcal{D}}
		\\ 
		= 
		\left\langle  f(\boldsymbol{y}), \boldsymbol{A}(\boldsymbol{y}):\mathrm{D} \boldsymbol{\psi} \right\rangle_{ \mathcal{D}} 
		~~
		\forall (\varphi, \boldsymbol{\psi}) \in 
		\mathbb{U}   \times \mathbb{G}.
	\end{multline}
	In accordance with (\ref{eq:Euler-Lagrange-stochastic-discrete-local}), we define the symmetric bilinear form $a_h: \left( L \mathbb{U}  \times L \mathbb{G} \right)^2  \rightarrow \mathbb{R}$ by
	%	\begin{equation*}
		%		a_h: \left( L \mathbb{U}  \times L \mathbb{G} \right)^2  \rightarrow \mathbb{R}
		%	\end{equation*}
	%	by
	\begin{multline}
		\label{def:stochastic-bilinear-a-discerete-local}
		a_h (\varphi, \boldsymbol{\psi}; \varphi^{\prime}, \boldsymbol{\psi}^{\prime}):=
		\left\langle \nabla\varphi - \boldsymbol{\psi} , \nabla\varphi^{\prime} - \boldsymbol{\psi}^{\prime}  \right\rangle_{\rho, \mathcal{D}}  
		+
		\left\langle  \nabla \times \boldsymbol{\psi}, \nabla \times \boldsymbol{\psi}^{\prime} \right\rangle_{\rho, \mathcal{D}}  
		\\
		+
		\left\langle  \boldsymbol{\psi}_{\boldsymbol{t}} , \boldsymbol{\psi}^{\prime}_{\boldsymbol{t}} \right\rangle _{\rho,\partial \mathcal{D}}
		+
		\sum_{ F \in \mathcal{F}^b} h_{K,F}^{-1}\left\langle  \boldsymbol{\psi}_{\boldsymbol{t}} , \boldsymbol{\psi}^{\prime}_{\boldsymbol{t}} \right\rangle _{\rho,F}
		+
		\left\langle  \boldsymbol{A}: \mathrm{D} \boldsymbol{\psi}, \boldsymbol{A}: \mathrm{D} \boldsymbol{\psi}^{\prime} \right\rangle_{\rho,\mathcal{D}}. 
	\end{multline}
	
	\vspace*{1mm}
	\begin{lemma}[discrete trace inequality]
		\label{lem:inverse-on-boundary}
		There exists $C_{\ref{ineq:inverse-on-boundary-local}}>0$ such that for every $\boldsymbol{\psi} \in \mathbb{G}$ and each face $F \in \mathcal{F}^b$ of an element $K$, the following holds
%		\begin{equation}
%			{\color{red}
%			\label{ineq:inverse-on-boundary-local-1}
%			\lVert \nabla_T (\boldsymbol{\psi} \cdot \boldsymbol{n}) \rVert_{L^2(F)} \leq C_{\ref{ineq:inverse-on-boundary-local-1}} h_K^{-\frac{1}{2}}
%			\lVert \mathrm{D} \boldsymbol{\psi} \rVert_{L^2( K)},
%		}
%		\end{equation}
		%and consequently, we have
		\begin{equation}
			\label{ineq:inverse-on-boundary-local}
			\lVert \left( \nabla (\boldsymbol{\psi} \cdot \boldsymbol{n}) \right)_{\boldsymbol{t}}  \rVert_{L^2(\partial \mathcal{D})} 
			\leq
			C_{\ref{ineq:inverse-on-boundary-local}}
			\sum_{F \in \mathcal{F}^b}
			h_{K,F}^{-\frac{1}{2}}
			\lVert \mathrm{D} \boldsymbol{\psi} \rVert_{L^2( K)} ,
		\end{equation}
	\end{lemma}
	
	\vspace*{1mm}
	\begin{proof}
		The claim follows directly from \cite[Lemma~12.8]{Ern-Guermond}.
	\end{proof}
	
	\vspace*{1mm}
	\begin{remark}[challenges with general functional spaces]
		\label{rem:discrete-vs.-countiuous-functional-space}
		%	The reason for considering (\ref{fun:cost-functional-discerete-local})–(\ref{def:stochastic-bilinear-a-discerete-local}) in $\left( L^2(\boldsymbol{\Gamma}) \otimes \mathbb{U} \right)  \times \left( L^2(\boldsymbol{\Gamma}) \otimes \mathbb{G} \right)$ instead of the more general space $\left( L^2(\boldsymbol{\Gamma}) \otimes H^1_0(D) \right)  \times \left( L^2(\boldsymbol{\Gamma}) \otimes H^1(D;\mathbb{R}^d) \right)$ is that the estimates (\ref{ineq:inverse-on-boundary-local-1}) and (\ref{ineq:inverse-on-boundary-local}) are not necessarily valid in the latter space. 
		We note that the formulation presented in (\ref{fun:cost-functional-discerete-local})-(\ref{eq:Euler-Lagrange-stochastic-discrete-local}) and (\ref{def:stochastic-bilinear-a-discerete-local}) could also be considered in the more general space $LH \times L \boldsymbol{H} $ instead of $L \mathbb{U} \times L \mathbb{G}$. 
		However, this level of generality cannot be pursued here because the estimate (\ref{ineq:inverse-on-boundary-local}) is not necessarily valid for every $\boldsymbol{\psi} \in \boldsymbol{H}$.
		%$\left( L^2(\boldsymbol{\Gamma}) \otimes H^1_0(D) \right)  \times \left( L^2(\boldsymbol{\Gamma}) \otimes H^1(D;\mathbb{R}^d) \right)$.	
	\end{remark}
	
	\vspace*{1mm}
	\begin{lemma}[discrete Miranda-Talenti type inequality]
		\label{lem:discrete-generalized-Maxwell}
		There exists a constant $C_{\ref{ineq:discrete-generalized-Maxwell-local}}>0$ such that for any $\boldsymbol{\psi} \in L \mathbb{G}$,
		\begin{multline}
			\label{ineq:discrete-generalized-Maxwell-local}
			\left( 1 - \frac{\varepsilon}{2} \right) 
			\left\| \mathrm{D} \boldsymbol{\psi} \right\|^2_{L^2_\rho(\boldsymbol{\Gamma}) \otimes L^2(\mathcal{D})}
			\\
			\leq
			\left\| \nabla \times \boldsymbol{\psi} \right\|^2_{L^2_\rho(\boldsymbol{\Gamma}) \otimes L^2(\mathcal{D})} 
			+
			C_{\ref{ineq:discrete-generalized-Maxwell-local}}
			\sum_{F \in \mathcal{F}^b} (1+ h_{K,F}^{-1})
			\lVert \boldsymbol{\psi}_{\boldsymbol{t}} \rVert^2_{ L^2_\rho(\boldsymbol{\Gamma}) \otimes L^2(F)}
			+
			\left\| \nabla \cdot \boldsymbol{\psi} \right\|^2_{ L^2_\rho(\boldsymbol{\Gamma}) \otimes L^2(\mathcal{\mathcal{D}})}.
		\end{multline}
		We recall that the parameter $\varepsilon$ is specified by the Cordes condition~(\ref{con:cordes}).
	\end{lemma}
	
	\vspace*{1mm}
	\begin{proof}
		The convexity of the domain $\mathcal{D}$ implies that, in Theorem~3.1.1.1 of \cite{Grisvard},
		$\left\langle (\text{tr }\mathscr{B}) \boldsymbol{\psi} \cdot \boldsymbol{n} , \boldsymbol{\psi} \cdot \boldsymbol{n} \right\rangle _{\partial \mathcal{D}} \leq 0$. Consequently, we have
		\begin{align*}
			\| \nabla \times & \boldsymbol{\psi} \|^2_{L^2_{\rho}(\boldsymbol{\Gamma}) \otimes L^2(\mathcal{D})} 
			+
			\left\| \nabla \cdot \boldsymbol{\psi} \right\|^2_{L^2_{\rho}(\boldsymbol{\Gamma}) \otimes L^2(\mathcal{D})}
			+
			C
			\lVert \boldsymbol{\psi}_{\boldsymbol{t}} \rVert^2_{L^2_{\rho}(\boldsymbol{\Gamma}) \otimes L^2(\partial \mathcal{D})}
			\\
			&
			\hspace*{-2mm}
			\geq
			\left\| \mathrm{D} \boldsymbol{\psi} \right\|^2_{L^2_{\rho}(\boldsymbol{\Gamma}) \otimes L^2(\mathcal{D})}
			-2 \sum_{F \in \mathcal{F}^b} \left\langle  \boldsymbol{\psi}_{\boldsymbol{t}} , \left( \nabla(\boldsymbol{\psi} \cdot \boldsymbol{n})\right)_{\boldsymbol{t}}  \right\rangle _{\rho, F}
			\\
			&
			\hspace*{-2mm}
			\geq
			\left\| \mathrm{D} \boldsymbol{\psi} \right\|^2_{L^2_{\rho}(\boldsymbol{\Gamma}) \otimes L^2(\mathcal{D})}
			-2 \sum_{F \in \mathcal{F}^b} \lVert \boldsymbol{\psi}_{\boldsymbol{t}} \rVert_{ L^2_{\rho}(\boldsymbol{\Gamma}) \otimes L^2(F)}
			\lVert \left(  \nabla (\boldsymbol{\psi} \cdot \boldsymbol{n}) \right)_{\boldsymbol{t}}  \rVert_{ L^2_{\rho}(\boldsymbol{\Gamma}) \otimes L^2(F)}
			\\
			&
			\hspace*{-2mm}
			\geq
			\left\| \mathrm{D} \boldsymbol{\psi} \right\|^2_{L^2_{\rho}(\boldsymbol{\Gamma}) \otimes L^2(\mathcal{D})}
			-2 \sum_{F \in \mathcal{F}^b} C_{\ref{ineq:inverse-on-boundary-local}} h_{K,F}^{-\frac{1}{2}} \lVert \boldsymbol{\psi}_{\boldsymbol{t}} \rVert_{ L^2_{\rho}(\boldsymbol{\Gamma}) \otimes L^2(F)}
			\lVert \mathrm{D} \boldsymbol{\psi} \rVert_{L^2_{\rho}(\boldsymbol{\Gamma}) \otimes  L^2( \mathcal{D})}
			\\
			&
			\hspace*{-2mm}
			\geq
			\left\| \mathrm{D} \boldsymbol{\psi} \right\|^2_{L^2_{\rho}(\boldsymbol{\Gamma}) \otimes L^2(\mathcal{D})}
			- 
			\frac{2}{\varepsilon}C_{\ref{ineq:inverse-on-boundary-local}}^2 \sum_{F \in \mathcal{F}^b} h_{K,F}^{-1} \lVert \boldsymbol{\psi}_{\boldsymbol{t}} \rVert^2_{ L^2_{\rho}(\boldsymbol{\Gamma}) \otimes L^2(F)} 
			- \frac{\varepsilon}{2} \left\| \mathrm{D} \boldsymbol{\psi} \right\|^2_{L^2_{\rho}(\boldsymbol{\Gamma}) \otimes L^2(\mathcal{D})},
		\end{align*}
		where the first inequality follows from Theorem~3.1.1.1 in \cite{Grisvard} for some positive constant $C$, the second from the Cauchy–Schwarz inequality, the third from (\ref{ineq:inverse-on-boundary-local}), and the last from the Young's inequality.
	\end{proof}
	
	\vspace*{1mm}
	\begin{lemma}
		\label{lem:discrete-Miranda-Talenti}
		Under Assumption~\ref{assum:A-f}~\textit{\ref{assum:matrix-A}}, there exists $C_{\ref{ineq:discrete-Miranda-Talenti-local}}>0$ such that 
		for any $\boldsymbol{\psi} \in L \mathbb{G}$,
		\begin{multline}
			\label{ineq:discrete-Miranda-Talenti-local}
			\left\| \nabla \times \boldsymbol{\psi} \right\|^2_{L^2_{\rho}(\boldsymbol{\Gamma}) \otimes L^2(\mathcal{D})} 
			+
			\sum_{F \in \mathcal{F}^b}
			(1+ h_{K,F}^{-1})
			\lVert \boldsymbol{\psi}_{\boldsymbol{t}} \rVert^2_{L^2_{\rho}(\boldsymbol{\Gamma}) \otimes L^2(F)}
			+
			\lVert \boldsymbol{A}:\mathrm{D} \boldsymbol{\psi} \rVert^2_{ L^2_{\rho}(\boldsymbol{\Gamma}) \otimes L^2(\mathcal{D})}
			\\
			\geq
			C_{\ref{ineq:discrete-Miranda-Talenti-local}}
			\lVert \mathrm{D} \boldsymbol{\psi} \rVert^2_{L^2_{\rho}(\boldsymbol{\Gamma}) \otimes L^2(\mathcal{D})} .
		\end{multline}
	\end{lemma}
	
	\vspace*{1mm}
	\begin{proof}
		The scaling function  $\gamma= \frac{\text{tr} \boldsymbol{A}}{\lvert \boldsymbol{A} \rvert^2}$ and the Cordes condition (\ref{con:cordes}) imply that
		\begin{equation}
			\label{ineq:discrete-gammaA-I}
			\lvert \gamma \boldsymbol{A} - \boldsymbol{I} \rvert^2 =
			d - \dfrac{\lvert \boldsymbol{A} \rvert^2}{\left( \text{tr} \boldsymbol{A} \right)^2 }
			\leq
			1-\varepsilon ,
			~ \text{for a.e. } (\boldsymbol{y}, \boldsymbol{x}) \in \boldsymbol{\Gamma} \times  \mathcal{D}. 
		\end{equation}
		Therefore, we have
		\begin{equation}
			\label{ineq:discrete-gammaA-Laplac}
			\lVert \left( \gamma \boldsymbol{A} - \boldsymbol{I} \right) : \mathrm{D} \boldsymbol{\psi} \rVert_{L^2_{\rho}(\boldsymbol{\Gamma}) \otimes L^2(\mathcal{D})}
			\leq \sqrt{1- \varepsilon} 
			\lVert \mathrm{D} \boldsymbol{\psi} \rVert_{L^2_{\rho}(\boldsymbol{\Gamma}) \otimes L^2(\mathcal{D})} .
		\end{equation}
		By adding and subtracting $\boldsymbol{I}: \mathrm{D} \boldsymbol{\psi}$, we have
		\begin{multline}
			\label{ineq:discrete-l-h-s-local}
			\left\| \nabla \times \boldsymbol{\psi} \right\|^2_{ L^2_{\rho}(\boldsymbol{\Gamma}) \otimes L^2(\mathcal{D})} 
			+
			C_{\ref{ineq:discrete-generalized-Maxwell-local}}
			\sum_{F \in \mathcal{F}^b}
			(1+ h_{K,F}^{-1})
			\lVert \boldsymbol{\psi}_{\boldsymbol{t}} \rVert^2_{L^2_{\rho}(\boldsymbol{\Gamma}) \otimes L^2(F)}
			+
			\lVert \gamma \boldsymbol{A}:\mathrm{D} \boldsymbol{\psi} \rVert^2_{L^2_{\rho}(\boldsymbol{\Gamma}) \otimes L^2(\mathcal{D})}
			\\
			\geq
			\bigg( 
				\Big(			
				\left\| \nabla \times \boldsymbol{\psi} \right\|^2_{L^2_{\rho}(\boldsymbol{\Gamma}) \otimes L^2(\mathcal{D})} 
				+
				C_{\ref{ineq:discrete-generalized-Maxwell-local}}
				\sum_{F \in \mathcal{F}^b}
				(1+ h_{K,F}^{-1}) 
				\lVert \boldsymbol{\psi}_{\boldsymbol{t}} \rVert^2_{ L^2_{\rho}(\boldsymbol{\Gamma}) \otimes L^2(F)}
				\\
				+
				\left\| \nabla \cdot \boldsymbol{\psi} \right\|^2_{ L^2_{\rho}(\boldsymbol{\Gamma}) \otimes L^2(\mathcal{D})}
			\Big)^{1/2} 
			%\\	
			-
			\lVert \left( \gamma \boldsymbol{A} - \boldsymbol{I} \right) : \mathrm{D} \boldsymbol{\psi} \rVert_{L^2_{\rho}(\boldsymbol{\Gamma}) \otimes L^2(\mathcal{D})}
			\bigg)^2
			\\
			\geq
			\left( \sqrt{1-\frac{\varepsilon}{2}} - \sqrt{1- \varepsilon}\right)^2 
			\left\| \mathrm{D} \boldsymbol{\psi} \right\|^2_{L^2_{\rho}(\boldsymbol{\Gamma}) \otimes L^2(\mathcal{D})}.
		\end{multline}
	\end{proof}
	
	\vspace*{1mm}
	\begin{lemma}[coercivity of $a_h(\cdot ; \cdot)$]
		\label{lem:coersivity-a_h}
		The stochastic bilinear form $a_h$ is coercive; i.e., there exists $C_{\ref{ineq:coercivity-a_h-local}}>0$, independent of $h$, such that for 
		any $(\varphi, \boldsymbol{\psi}) \in  L \mathbb{U}  \times L \mathbb{G} $ the following holds
		\begin{multline}
			\label{ineq:coercivity-a_h-local}
			a_h (\varphi, \boldsymbol{\psi}; \varphi, \boldsymbol{\psi}) := 
			\lVert \nabla \varphi - \boldsymbol{\psi} \rVert^2_{L^2_{\rho}(\boldsymbol{\Gamma}) \otimes L^2(\mathcal{D})}
			+
			\lVert \nabla \times \boldsymbol{\psi} \rVert^2_{L^2_{\rho}(\boldsymbol{\Gamma}) \otimes L^2(\mathcal{D})}
			\\ 
			+
			\sum_{F \in \mathcal{F}^b}
			\left( 1 + h_{K,F}^{-1} \right) 
			\lVert \boldsymbol{\psi}_{\boldsymbol{t}} \rVert^2_{L^2_{\rho}(\boldsymbol{\Gamma}) \otimes L^2(F)}
			+
			\lVert \boldsymbol{A}: \mathrm{D} \boldsymbol{\psi} \rVert_{L^2_{\rho}(\boldsymbol{\Gamma}) \otimes L^2(\mathcal{D})}^2
			\\
			\geq
			C_{\ref{ineq:coercivity-a_h-local}}
			\left\| (\varphi, \boldsymbol{\psi})  \right\| ^2_{L^2_{\rho}(\boldsymbol{\Gamma}) \otimes H^1_h(\mathcal{D})}.
		\end{multline} 
	\end{lemma}
	
	\vspace*{1mm}
	\begin{proof}
		Lemma~\ref{lem:discrete-Miranda-Talenti}, together with Young's inequality and Poincaré's inequality, implies that
		\begin{align*}
			\frac{C_P C_{\ref{ineq:discrete-Miranda-Talenti-local}}}{2}& \lVert \nabla \varphi - \boldsymbol{\psi} \rVert^2_{L^2_{\rho}(\boldsymbol{\Gamma}) \otimes L^2(\mathcal{D})}
			+
			\lVert \nabla \times \boldsymbol{\psi} \rVert^2_{L^2_{\rho}(\boldsymbol{\Gamma}) \otimes L^2(\mathcal{D})} 
			\\
			+
			&(1+ C_{\ref{ineq:discrete-Miranda-Talenti-local}})
			\sum_{F \in \mathcal{F}^b}
			(1 + h_{K,F}^{-1})\lVert \boldsymbol{\psi}_{\boldsymbol{t}} \rVert^2_{L^2_{\rho}(\boldsymbol{\Gamma}) \otimes L^2(F)}
			+
			\lVert \boldsymbol{A} : \mathrm{D} \boldsymbol{\psi} \rVert_{L^2_{\rho}(\boldsymbol{\Gamma}) \otimes L^2(\mathcal{D})}^2
			\\
			&{\hspace*{-5.5mm}}\geq
			\frac{C_P C_{\ref{ineq:discrete-Miranda-Talenti-local}}}{2} 
			\lVert \nabla \varphi - \boldsymbol{\psi} \rVert^2_{L^2_{\rho}(\boldsymbol{\Gamma}) \otimes L^2(\mathcal{D})}
			+
			C_{\ref{ineq:discrete-Miranda-Talenti-local}} 
			\lVert \boldsymbol{\psi} \rVert^2_{L^2_{\rho}(\boldsymbol{\Gamma}) \otimes H^1_h(\mathcal{D})} 
			\\
			&{\hspace*{-5.5mm}}\geq
			\frac{C_P C_{\ref{ineq:discrete-Miranda-Talenti-local}} }{4} 
			\lVert \nabla \varphi \rVert_{L^2_{\rho}(\boldsymbol{\Gamma}) \otimes L^2(\mathcal{D})}^2 
			-
			\frac{C_P C_{\ref{ineq:discrete-Miranda-Talenti-local}}}{2} \lVert \boldsymbol{\psi} \rVert_{L^2_{\rho}(\boldsymbol{\Gamma}) \otimes L^2(\mathcal{D})}^2
			+
			C_{\ref{ineq:discrete-Miranda-Talenti-local}} 
			\lVert \boldsymbol{\psi} \rVert^2_{L^2_{\rho}(\boldsymbol{\Gamma}) \otimes H^1_h(\mathcal{D})}
			\\
			&{\hspace*{-5.5mm}} \geq
			\frac{C_P C_{\ref{ineq:discrete-Miranda-Talenti-local}} }{4} 
			\lVert \nabla \varphi \rVert_{L^2_{\rho}(\boldsymbol{\Gamma}) \otimes L^2(\mathcal{D})}^2 
			+
			\frac{C_{\ref{ineq:discrete-Miranda-Talenti-local}}}{2} 
			\lVert \boldsymbol{\psi} \rVert^2_{L^2_{\rho}(\boldsymbol{\Gamma}) \otimes H^1_h(\mathcal{D})}
			,
		\end{align*} 
		where $C_P>0$ is the constant of Poincaré's inequality:
		\begin{equation*}
			\lVert \mathrm{D} \boldsymbol{\psi} \rVert^2_{L^2_\rho(\boldsymbol{\Gamma}) \otimes L^2(\mathcal{D})}
			+
			\lVert \boldsymbol{\psi}_{\boldsymbol{t}} \rVert^2_{ L^2_\rho(\boldsymbol{\Gamma}) \otimes L^2(\partial \mathcal{D})}
			\geq
			C_P
			\lVert \boldsymbol{\psi} \rVert^2_{L^2_\rho(\boldsymbol{\Gamma}) \otimes L^2(\mathcal{D})}.
		\end{equation*}
	\end{proof}
	
	\vspace*{1mm}
	%{\color{red}
		%\begin{remark}[{{\bf coercivity of $a_h(\cdot, \cdot)$}} \bf on $u+ {\color{red}W}$]
		%	\label{rem:coercivity-with-solution}
		%	As mentioned in Remark \ref{rem:discrete-vs.-countiuous-functional-space}, it is possible to define the bilinear form $a_h$ on the more general space $\left( L^2(\boldsymbol{\Gamma}) \otimes H^1_0(D) \right)  \times \left( L^2(\boldsymbol{\Gamma}) \otimes H^1(D;\mathbb{R}^d) \right)$ with the same formulation. However, for the strong solution $u$ of (\ref{eq:nondivergence-random-variable}), the vanishing of the challenging term $\left( \boldsymbol{t}_{D} \nabla u\right) \vert_{\partial D}=\boldsymbol{0}$ simplifies the problem. Consequently, by Lemma \ref{lem:coersivity-a_h}, the bilinear form $a_h$ remains coercive on $u + {\color{red}W}$. Specifically, for any $(\varphi, \boldsymbol{\psi}) \in  \left( L^2(\boldsymbol{\Gamma}) \otimes \mathbb{U} \right)  \times \left(  L^2(\boldsymbol{\Gamma}) \otimes \mathbb{G}\right) $, the following holds
		%	\begin{equation}
			%		\label{ineq:coercivity-with-solution}
			%			a_h (u-\varphi, \nabla u - \boldsymbol{\psi}; u - \varphi, \nabla u - \boldsymbol{\psi})
			%			\geq
			%			C_{\ref{ineq:coercivity-a_h-local}}
			%			\left\| (u - \varphi, \nabla u - \boldsymbol{\psi})  \right\| ^2_{L^2(\boldsymbol{\Gamma}) \otimes H^1_h(D)}.
			%	\end{equation}
		%\end{remark}}
		The coercivity condition (\ref{ineq:coercivity-a_h-local}) alone does not guarantee the well-posedness of problem (\ref{eq:Euler-Lagrange-stochastic-discrete-local}). This is because the solution space $ L \mathbb{U} \times L \mathbb{G}$ is not finite-dimensional; while the deterministic component of the space has finite dimension, the stochastic component remains infinite-dimensional. Therefore, to ensure the well-posedness of the solution, the continuity of the bilinear form $a_h$ is also required.
		
		\vspace*{1mm}
		\begin{remark}[continuity of $a_h(\cdot ; \cdot)$]
			It is straightforward to demonstrate that the mesh-dependent bilinear form	$a_h(\cdot; \cdot)$ is continuous on $L H \times L \boldsymbol{H}$
			under the mesh-dependent norm defined in (\ref{eq:discrete-norm-U-times-G}); 
			i.e., there exist $C_{\ref{ineq:continuity-a_h}}> 0$ such that for 
			any $(\varphi, \boldsymbol{\psi}), (\varphi^{\prime}, \boldsymbol{\psi}^{\prime}) \in L H \times L \boldsymbol{H} $ the following holds
			\begin{equation}
				\label{ineq:continuity-a_h}
				\left|  a_h (\varphi, \boldsymbol{\psi}; \varphi^\prime, \boldsymbol{\psi})^\prime \right| 
				\leq
				C_{\ref{ineq:continuity-a_h}}
				\left\|  (\varphi, \boldsymbol{\psi}) \right\|_{L^2_{\rho}(\boldsymbol{\Gamma}) \otimes
					H^1_h(\mathcal{D})   
				}
				\left\|  (\varphi^{\prime}, \boldsymbol{\psi}^\prime) \right\|_{ L^2_{}\rho(\boldsymbol{\Gamma}) \otimes H^1_h(\mathcal{D})  
				}.
			\end{equation}
			However, we emphasize that this continuity holds only with respect to the mesh-dependent norm and not with the continuous norm.
		\end{remark}
		
		\vspace*{1mm}
		\begin{remark}[well-posedness of (\ref{eq:Euler-Lagrange-stochastic-discrete-local})]
			The coercivity~(\ref{ineq:coercivity-a_h-local}), combined with the continuity~(\ref{ineq:continuity-a_h}), ensures the well-posedness of problem~(\ref{eq:Euler-Lagrange-stochastic-discrete-local}). Moreover, using similar arguments, the well-posedness of problem~(\ref{eq:Euler-Lagrange-stochastic-discrete-local-parametric}) can also be established.
		\end{remark}
		
		\vspace*{1mm}
			For any $(\varphi, \boldsymbol{\psi}) \in L H \times L \boldsymbol{H}$, let $\mathcal{I}_{\mathbb{U}} \varphi$ and $\mathcal{I}_{\mathbb{G}} \boldsymbol{\psi}$ represent the nodal interpolations of $\varphi$ in $\mathbb{U}$ and $\boldsymbol{\psi}$ in $\mathbb{G}$ respectively.
			
			\vspace*{1mm}
			\begin{lemma}[interpolation error estimate {\cite[Theorem~3.1.6]{Ciarlet}}]
				\label{lem:interpolation-error}
				For any $(\varphi, \boldsymbol{\psi}) \in \left( L^2_{\rho}(\boldsymbol{\Gamma}) \otimes H^s (\mathcal{D}) \right) \times \left( L^2_{\rho}(\boldsymbol{\Gamma}) \otimes H^s (\mathcal{D}; \mathbb{R}^d) \right) $ with $s>m$, there exists $C_{\ref{ineq:interpolation-error}}>0$ such that
				\begin{equation}
					\begin{split}
						\label{ineq:interpolation-error}
						\lVert \varphi - \mathcal{I}_{\mathbb{U}} \varphi \rVert_{L^2_{\rho}(\boldsymbol{\Gamma}) \otimes H^m(\mathcal{D})} 
						&\leq
						C_{\ref{ineq:interpolation-error}}h^{l - m} \lVert \varphi \rVert_{L^2_{\rho}(\boldsymbol{\Gamma}) \otimes H^s(\mathcal{D})},
						\\
						\lVert \boldsymbol{\psi} - \mathcal{I}_{\mathbb{G}} \boldsymbol{\psi} \rVert_{L^2_{\rho}(\boldsymbol{\Gamma}) \otimes H^m(\mathcal{D})} 
						&\leq
						C_{\ref{ineq:interpolation-error}}h^{l - m} \lVert \boldsymbol{\psi} \rVert_{L^2_{\rho}(\boldsymbol{\Gamma}) \otimes H^s(\mathcal{D})},
					\end{split}
				\end{equation}
				where $l:= \min\left\lbrace k+1 , s\right\rbrace $. We recall that $k$ denotes the degree of the polynomials in the finite element spaces defined in (\ref{def:Galerkin-spaces}).
			\end{lemma}
			
			\vspace*{1mm}
			\begin{remark}[interpolation error of strong solution]
				\label{rem:interpolation-error-discrete-norm}
				For the strong solution $u$ of (\ref{eq:nondivergence-random-variable}), the vanishing tangential trace of $\nabla u$, $\left( \left(  \nabla u \right)_{\boldsymbol{t}}  \right) =\boldsymbol{0}$, and consequently, the vanishing tangential trace of $ \mathcal{I}_{\mathbb{G}} \nabla u$, $\left( \left(  \mathcal{I}_{\mathbb{G}} \nabla u\right)_{\boldsymbol{t}} \right) =\boldsymbol{0}$, ensure that the interpolation error of $\nabla u$ with $L^2_{\rho}(\boldsymbol{\Gamma}) \otimes H^1(\mathcal{D})$ and $L^2_{\rho}(\boldsymbol{\Gamma}) \otimes H^1_h(\mathcal{D})$ norms are equivalent. Specifically, when $u$ satisfies $u \in L^2_{\rho}(\boldsymbol{\Gamma}) \otimes H^{2+ \varrho}(\mathcal{D})$ for some real $\varrho >0$, Lemma \ref{lem:interpolation-error} implies the existence of a constant $C_{\ref{ineq:interpolation-error-discrete-norm}}>0$ such that
				\begin{equation}
					\label{ineq:interpolation-error-discrete-norm}
					\lVert (u, \nabla u) - (\mathcal{I}_{\mathbb{U}}u , \mathcal{I}_{\mathbb{G}} \nabla u) \rVert_{L^2_{\rho}(\boldsymbol{\Gamma}) \otimes H^1_h(\mathcal{D})}
					\leq
					C_{\ref{ineq:interpolation-error-discrete-norm}}
					h^{\min\left\lbrace k, \varrho \right\rbrace } \lVert u \rVert_{L^2_{\rho}(\boldsymbol{\Gamma}) \otimes H^{2 + \varrho}(\mathcal{D})}.
				\end{equation}	
			\end{remark}

			\vspace*{1mm}
			\begin{theorem}[a priori error estimate]
				\label{lem:a-priori-error-estimate}
				Let $\mathcal{T}$ be in a collection $\mathfrak{T}$ of
				shape-regular conforming simplicial meshes on the polyhedral domain $\mathcal{D} \subset \mathbb{R}^d$. Suppose that the strong solution $u$ of (\ref{eq:nondivergence-random-variable}) satisfies $u \in L^2_{\rho}(\boldsymbol{\Gamma}) \otimes H^{2+ \varrho}(\mathcal{D})$ for some real $\varrho >0$, and $(u_{\mathbb{U}} , \boldsymbol{g}_{\mathbb{G}})$ is the unique solution of the semi-discrete problem (\ref{eq:Euler-Lagrange-stochastic-discrete-local}). Then
				\begin{equation}
					\label{ineq:a-priori-error-estimate}
					\lVert (u, \nabla u) - (u_{\mathbb{U}} , \boldsymbol{g}_{\mathbb{G}}) \rVert_{L^2_{\rho}(\boldsymbol{\Gamma}) \otimes H^1_h(\mathcal{D})}
					\leq
					C_{\ref{ineq:a-priori-error-estimate}}
					h^{\min\left\lbrace k, \varrho \right\rbrace } \lVert u \rVert_{L^2_{\rho}(\boldsymbol{\Gamma}) \otimes H^{2 + \varrho}(\mathcal{D})}.
				\end{equation}
			\end{theorem}
			
			\vspace*{1mm}
			\begin{proof}
				%		{ \color{red}Combining the results of Lemma \ref{lem:quasi-optimality-discrete-norm} and Lemma \ref{lem:interpolation-error-discrete-norm} demonstrates the claim.}
				%	
				The triangle inequality implies that
				\begin{multline}
					\label{ineq:triangle-discrete-norm}
					\lVert (u, \nabla u) - (u_{\mathbb{U}}, \boldsymbol{g}_{\mathbb{G}}) \rVert_{L^2_{\rho}(\boldsymbol{\Gamma}) \otimes H^1_h(\mathcal{D})}
					\\
					\leq  
					\lVert (u, \nabla u) - (\mathcal{I}_{\mathbb{U}}u, \mathcal{I}_{\mathbb{G}} \nabla u) \rVert_{L^2_{\rho}(\boldsymbol{\Gamma}) \otimes H^1_h(\mathcal{D})}
					+
					\lVert (\mathcal{I}_{\mathbb{U}}u, \mathcal{I}_{\mathbb{G}} \nabla u) - (u_{\mathbb{U}}, \boldsymbol{g}_{\mathbb{G}}) \rVert_{L^2_{\rho}(\boldsymbol{\Gamma}) \otimes H^1_h(\mathcal{D})}.
				\end{multline}
				We provide an upper bound for the first term on the right-hand side of (\ref{ineq:triangle-discrete-norm}) using the error estimate (\ref{ineq:interpolation-error-discrete-norm}), as stated in Remark~\ref{rem:interpolation-error-discrete-norm}. To bound the second term, we note that 
				\begin{equation*}
					a_h\left( (u, \nabla u)- (u_{\mathbb{U}} , \boldsymbol{g}_{\mathbb{G}}) ; (\varphi, \boldsymbol{\psi}) \right) = 0, \quad \forall  (\varphi, \boldsymbol{\psi}) \in  L \mathbb{U} \times L \mathbb{G} .
				\end{equation*}
				Therefore we have
				\begin{align*}
					\lVert (\mathcal{I}_{\mathbb{U}}&u, \mathcal{I}_{\mathbb{G}} \nabla u) - (u_{\mathbb{U}}, \boldsymbol{g}_{\mathbb{G}}) \rVert_{L^2_{\rho}(\boldsymbol{\Gamma}) \otimes H^1_h(\mathcal{D})}^2
					\\
					&{\hspace*{-3.5mm}}\leq
					C_{\ref{ineq:coercivity-a_h-local}}^{-1} 
					~
					a_h\left(  (\mathcal{I}_{\mathbb{U}}u, \mathcal{I}_{\mathbb{G}} \nabla u) - (u_{\mathbb{U}}, \boldsymbol{g}_{\mathbb{G}}); (\mathcal{I}_{\mathbb{U}}u, \mathcal{I}_{\mathbb{G}} \nabla u) - (u_{\mathbb{U}}, \boldsymbol{g}_{\mathbb{G}}) \right)
					\\
					&{\hspace*{-3.5mm}}=
					C_{\ref{ineq:coercivity-a_h-local}}^{-1}
					a_h\left(  (\mathcal{I}_{\mathbb{U}}u, \mathcal{I}_{\mathbb{G}} \nabla u) -
					(u, \nabla u) ; (\mathcal{I}_{\mathbb{U}}u, \mathcal{I}_{\mathbb{G}} \nabla u) - (u_{\mathbb{U}}, \boldsymbol{g}_{\mathbb{G}}) \right)
					\\
					&{\hspace*{-3.5mm}}\leq
					\dfrac{C_{\ref{ineq:continuity-a_h}}}{C_{\ref{ineq:coercivity-a_h-local}}} 
					\lVert (\mathcal{I}_{\mathbb{U}}u, \mathcal{I}_{\mathbb{G}} \nabla u) - (u, \nabla u) \rVert_{L^2_{\rho}(\boldsymbol{\Gamma}) \otimes H^1_h(\mathcal{D})}
					\lVert (\mathcal{I}_{\mathbb{U}}u, \mathcal{I}_{\mathbb{G}} \nabla u) - (u_{\mathbb{U}}, \boldsymbol{g}_{\mathbb{G}}) \rVert_{L^2_{\rho}(\boldsymbol{\Gamma}) \otimes H^1_h(\mathcal{D})},
				\end{align*}
				which implies that
				\begin{equation*}
					\label{ineq:near-optimal-discrete-norm}
					\lVert (\mathcal{I}_{\mathbb{U}}u, \mathcal{I}_{\mathbb{G}} \nabla u) - (u_{\mathbb{U}}, \boldsymbol{g}_{\mathbb{G}}) \rVert_{L^2_{\rho}(\boldsymbol{\Gamma}) \otimes H^1_h(\mathcal{D})}
					\leq
					\dfrac{C_{\ref{ineq:continuity-a_h}}}{C_{\ref{ineq:coercivity-a_h-local}}} 
					\lVert (\mathcal{I}_{\mathbb{U}}u, \mathcal{I}_{\mathbb{G}} \nabla u) - (u, \nabla u) \rVert_{L^2_{\rho}(\boldsymbol{\Gamma}) \otimes H^1_h(\mathcal{D})}.
				\end{equation*}
				Applying the error estimate (\ref{ineq:interpolation-error-discrete-norm}) once again completes the proof.
			\end{proof} 
			
			\vspace*{1mm}
			The content and derivations in the next two sections are adapted from \cite{Babuska-Nobile-Tempone}, with modifications and extensions where necessary to suit the current context.
	
	\vspace*{2mm}
	\subsection{Collocation method}
	\label{sec:Collocation-method}
	In this part, we review a well-established collocation method and adopt the notations from \cite{Babuska-Nobile-Tempone}.
	For certain classes of problems, when the solution of stochastic partial differential equation has a very smooth dependence on the input random variables, it is reasonable to use a global polynomial approximation in the parameter space $L^2_{\rho}(\boldsymbol{\Gamma})$. 
	Since in this work, the analytic dependence of the solution $u$ and consequently the gradient recovery $\boldsymbol{g}$ with respect to the random variables $y_n \in \Gamma_n,~ n=1, \cdots, N$, is fulfilled by Assumption~\ref{assum:convergence-A-f}, for the purpose of approximation, the stochastic domain $\boldsymbol{\Gamma}$ does not need to be refined. 
	For a global approximation, we define $\mathcal{P}_{\boldsymbol{p}}(\boldsymbol{\Gamma}) \subset L^2{\rho}(\boldsymbol{\Gamma})$ as the span of tensor product polynomials with degrees at most $\boldsymbol{p} = (p_1, \cdots, p_N)$. Specifically, $\mathcal{P}_{\boldsymbol{p}}(\boldsymbol{\Gamma})= \bigotimes_{n=1}^N\mathcal{P}_{p_n}(\Gamma_n)$, where 
	\begin{equation}
		\label{de-space:spectral}
		\mathcal{P}_{p_n}(\Gamma_n):=\text{span}(y_n^j, ~ j=0,\cdots,p_n), \quad n=1, \cdots, N.
	\end{equation}
	Hence, the dimension of $\mathcal{P}_{\boldsymbol{p}}$ is $N_{\boldsymbol{p}} = \prod_{n=1}^N(p_n + 1)$.
	
	By setting $\left( \mathcal{P}_{\boldsymbol{p}}(\boldsymbol{\Gamma}) \otimes \mathbb{U} \right)   \times \left( \mathcal{P}_{\boldsymbol{p}}(\boldsymbol{\Gamma}) \otimes  \mathbb{G} \right) \subset L \mathbb{U} \times L \mathbb{G}$ as the finite element space, the fully-discrete counterpart of (\ref{eq:Euler-Lagrange-stochastic-discrete-local}) consists in finding $(u_{\boldsymbol{p},\mathbb{U}}, \boldsymbol{g}_{\boldsymbol{p},\mathbb{G}}) \in \left(  \mathcal{P}_{\boldsymbol{p}}(\boldsymbol{\Gamma}) \otimes \mathbb{U} \right)   \times \left( \mathcal{P}_{\boldsymbol{p}}(\boldsymbol{\Gamma}) \otimes \mathbb{G} \right) $ such that
	\begin{multline}
		\label{eq:Euler-Lagrange-discrete}
		\left\langle \nabla u_{\boldsymbol{p},\mathbb{U}} - \boldsymbol{g}_{\boldsymbol{p},\mathbb{G}}, \nabla\varphi - \boldsymbol{\psi} \right\rangle_{\rho,\mathcal{D}}  
		+
		\left\langle  \nabla \times \boldsymbol{g}_{\boldsymbol{p},\mathbb{G}}, \nabla \times \boldsymbol{\psi} \right\rangle_{\rho, \mathcal{D}} 
		\\
		+
		\sum_{F \in \mathcal{F}^b} \left( 1 + h_{K,F}^{-1}\right) \left\langle   \left( \boldsymbol{g}_{\boldsymbol{p},\mathbb{G}}\right) _{\boldsymbol{t}} , \boldsymbol{\psi}_{\boldsymbol{t}} \right\rangle _{\rho, F}
		+
		\left\langle  \boldsymbol{A}: \mathrm{D} \boldsymbol{g}_{\boldsymbol{p},\mathbb{G}}, \boldsymbol{A}: \mathrm{D} \boldsymbol{\psi} \right\rangle_{\rho, \mathcal{D}} 
		\\
		= 
		\left\langle  f, \boldsymbol{A}: \mathrm{D} \boldsymbol{\psi} \right\rangle_{\rho, \mathcal{D}} 
		~~
		\forall (\varphi, \boldsymbol{\psi}) \in 
		\left( \mathcal{P}_{\boldsymbol{p}}(\boldsymbol{\Gamma}) \otimes \mathbb{U} \right)  \times \left(  \mathcal{P}_{\boldsymbol{p}}(\boldsymbol{\Gamma}) \otimes \mathbb{G}\right)  .
	\end{multline}
	In general, problem (\ref{eq:Euler-Lagrange-discrete}) results in a fully coupled system of linear equations with dimension $N_{h,k} \times N_{\boldsymbol{p}}$, where parallel computation can not be effectively employed and highly efficient strategies are lacking for its numerical solution.
	In contrast, the collocation method only requires solving $N_{\boldsymbol{p}}$ uncoupled linear systems, each of dimension $N_{h,k}$, making it fully parallelizable.
	To benefit from the computational advantages of decoupled systems for equations with more general data, we adopt the collocation method to estimate the solution within the stochastic domain.	
	
	By selecting a set $\left\lbrace \boldsymbol{y}_1 , \cdots, \boldsymbol{y}_{N_{\boldsymbol{p}}}\right\rbrace $ as collocation points in $\boldsymbol{\Gamma}$, the fully discrete approximation is expressed as
	\begin{equation}
		\label{def:discrete-solution}
		u_{ \boldsymbol{p} , \mathbb{U}}(\boldsymbol{y}, \boldsymbol{x}) := 
		\sum_{m=1}^{N_{\boldsymbol{p}}} u_{\mathbb{U}}(\boldsymbol{y}_m, \boldsymbol{x})
		\ell_m(\boldsymbol{y}), 
		\quad
		\boldsymbol{g}_{ \boldsymbol{p} , \mathbb{G}}(\boldsymbol{y}, \boldsymbol{x}) := 
		\sum_{m=1}^{N_{\boldsymbol{p}}} \boldsymbol{g}_{\mathbb{G}}(\boldsymbol{y}_m, \boldsymbol{x})
		\ell_m(\boldsymbol{y}),
	\end{equation}
	where $\ell_m$ are multivariate Lagrange basis functions associated with the interpolation points $\boldsymbol{y}_m$, and $\left( u_{\mathbb{U}}(\boldsymbol{y}_m, \cdot), \boldsymbol{g}_{\mathbb{G}}(\boldsymbol{y}_m, \cdot) \right) $ is the solution of (\ref{eq:Euler-Lagrange-stochastic-discrete-local-parametric}) with $\boldsymbol{y} = \boldsymbol{y}_m$. 
	For convenience, we introduce the Lagrange interpolation operator $\mathcal{I}_{\boldsymbol{p}}: C^0(\boldsymbol{\Gamma}; V) \rightarrow \mathcal{P}_{\boldsymbol{p}}(\boldsymbol{\Gamma}) \otimes V$ defined for any $v \in C^0(\boldsymbol{\Gamma}; V)$ as
	\begin{equation}
		\label{def:interpolation-op}
		\mathcal{I}_{\boldsymbol{p}} v(\boldsymbol{y}) := \sum_{m=1}^{N_{\boldsymbol{p}}}{v(\boldsymbol{y}_m) \ell_m(\boldsymbol{y})}.
	\end{equation}
	Consequently, we have $u_{\boldsymbol{p}, \mathbb{U}} = \mathcal{I}_{\boldsymbol{p}}{u_\mathbb{U}}$ and $ \boldsymbol{g}_{\boldsymbol{p} , \mathbb{G}} = \mathcal{I}_{\boldsymbol{p}} \boldsymbol{g_\mathbb{G}}$.
	
	To determine how the collocation points are set, we first introduce an auxiliary probability density function $\hat{\rho}:\boldsymbol{\Gamma} \rightarrow \mathbb{R}^{+}$, defined as
	\begin{equation*}
		\label{def:join-rho-idep}
		\hat{\rho}(\boldsymbol{y}):= \prod_{n=1}^N \hat{\rho}_n(y_n), ~~ \forall \boldsymbol{y} \in \boldsymbol{\Gamma}, 
	\end{equation*}
	such that
	\begin{equation}
		\label{def:join-rho}
		\left\|  \dfrac{\rho}{\hat{\rho}} 
		\right\| _{L^{\infty}(\boldsymbol{\Gamma})} < \infty.
	\end{equation} 
	This $\hat{\rho}$ represents the joint probability density of $N$ independent random variables. When the random variables $\left\lbrace y_n\right\rbrace _{n=1}^N$ are non-independent, $\rho(\boldsymbol{y})$ and $\hat{\rho}(\boldsymbol{y})$ are distinct functions. However, in the case of independent random variables, these two functions may be equal.
	
	We will select the collocation points as the zeros of suitable orthogonal polynomials. Specifically, for each dimension $n= 1, \cdots, N$, let $y_{n, m_{n}}$, $1 \leq m_n \leq p_n+1$, represent the the $p_n+1$ roots of the orthogonal polynomial $q_{p_n+1}$ with respect to the weight $\hat{\rho}_n$. This polynomial satisfies the condition
	\begin{equation*}
		\int_{\Gamma_n} q_{p_n+1}(y) v(y) \hat{\rho}_n(y) dy = 0 ~~
		\forall v \in \mathcal{P}_{p_{n}}(\Gamma_n).
	\end{equation*}	
	We also introduce the Lagrange basis $\left\lbrace \ell_{n,m_n}\right\rbrace_{m_n=1}^{p_n+1} \subset \mathcal{P}_{p_n}(\Gamma_n)$ as
	\begin{equation*}
		\ell_{n,m_n}(y_n):=\prod_{\stackbin[j\neq m_n]{}{j=1}}^{p_n+1}\frac{(y_n - y_{n,j})}{(y_{n,m_n} - y_{n,j})},
	\end{equation*}
	and the weights of the Gaussian quadrature formula $w_{n, m_n}$ as
	\begin{equation*}
		w_{n,m_n}:= \int_{\Gamma_n} \ell_{n,m_n}(y) \hat{\rho}_n(y)dy.
	\end{equation*}
	
	\vspace*{1mm}
	\begin{remark}
		\label{rem:equivalence-weights}
		For each $n=1,\dots, N$, the property $\sum_{j= 1}^{p_n+1}\ell_{n,j}(y) = 1$, implies that
		\begin{equation*}
			w_{n,m_n}= \int_{\Gamma_n} \ell_{n,m_n}(y) \hat{\rho}_n(y) dy = 
			\int_{\Gamma_n} \ell_{n,m_n}(y) 
			\left( \sum_{j= 1}^{p_n+1}\ell_{n,j}(y) \right)
			\hat{\rho}_n(y) dy.
		\end{equation*}
		Using the orthogonality condition,  $\int_{\Gamma_n}\ell_{n,m_n}(y) \ell_{n,j}(y) \hat{\rho}_n(y) dy =  0 ~ \forall j \neq  m_n$ we obtain
		\begin{equation}
			\label{eq:weights}
			w_{n,m_n}= \int_{\Gamma_n} \ell_{n,m_n}(y) \hat{\rho}_n(y) dy = 
			\int_{\Gamma_n} \ell_{n,m_n}^2(y) \hat{\rho}_n(y) dy.
		\end{equation}
	\end{remark}
	
	\vspace*{1mm}
	Standard choices for $\hat{\rho}_n$, such as constant or Gaussian distributions, lead to well-known roots of the polynomial $q_{p_{n+1}}$ and corresponding weights $w_{n, m_n}$, which are tabulated with full accuracy and do not require computation.
	
	For any index vector $( m_1, \cdots , m_N ) $ we associate a global index $m$ defined as
	\begin{equation*}
		m = m_1 + \sum_{i=1}^{N-1} (m_{i+1}-1) \prod_{j=1}^i (p_{j}+1).
	\end{equation*}	
	%	\begin{equation*}
		%		m= m_1 + (p_1+1)(m_2-1) + (p_1 +1) (p_2 + 1) (m_3 - 1) + \cdots 
		%		+
		%		\prod_{i=1}^{N-1}(p_i +1)(m_N-1). 
		%	\end{equation*}
	We then denote the corresponding collocation point by
	\begin{equation*}
		\boldsymbol{y}_m = ( y_{1,m_1}, y_{2,m_2} , \cdots , y_{N,m_N} ) \in \boldsymbol{\Gamma},
	\end{equation*}
	and set \begin{equation*}
		\ell_m(\boldsymbol{y}) := \prod_{n=1}^N \ell_{n, m_n} (y_n), 
		\quad 
		w_m := \prod_{n=1}^{N}w_{n, m_n}.
	\end{equation*}
	For any continuous function $\boldsymbol{\phi} : \boldsymbol{\Gamma} \rightarrow \mathbb{R}^k, ~ k\in \mathbb{N}$, the mean value is defined as $\bar{\boldsymbol{\phi}}:= \mathrm{E}_{\hat{\rho}} \left[  \frac{\rho}{\hat{\rho}} \boldsymbol{\phi} \right]$.
	When $\rho / \hat{\rho}$ is a smooth function, the Gauss quadrature formula can be used to approximate the integral, resulting in the following expressions:
	\begin{equation}
		\label{def:mean-value-finite-element}
		\bar{u}_{\boldsymbol{p} ,\mathbb{U}}(\boldsymbol{x}) \approx \sum_{k=1}^{N_{\boldsymbol{p}}} \frac{\rho(\boldsymbol{y}_k)}{\hat{\rho}(\boldsymbol{y}_k)} u_{\mathbb{U}}(\boldsymbol{y}_k, \boldsymbol{x}) w_k, 
		\quad
		\bar{\boldsymbol{g}}_{\boldsymbol{p} ,\mathbb{G}}(\boldsymbol{x}) \approx 
		\sum_{k=1}^{N_{\boldsymbol{p}}}
		\frac{\rho(\boldsymbol{y}_k)}{\hat{\rho}(\boldsymbol{y}_k)}
		\boldsymbol{g}_{\mathbb{G}}(\boldsymbol{y}_k, \boldsymbol{x}) w_k.
	\end{equation}
	However, if $\rho / \hat{\rho}$ exhibits discontinuities or singularities, the Gauss quadrature formula may no longer suffice. In such cases, a more appropriate quadrature formula must be employed to accurately compute  $\bar{u}_{\boldsymbol{p} ,\mathbb{U}}$ and $\bar{\boldsymbol{g}}_{\boldsymbol{p},\mathbb{G}}$, ensuring the effective handling of these irregularities.	
	
	\vspace*{2mm}
	\section{A priori error estimate and convergence}
	\label{sec:A-priori-error-convergence}
	
	To establish an a priori estimate for the approximation error, the total error is divided into two components: the error of finite element approximation in the deterministic domain, which can be estimated by (\ref{ineq:a-priori-error-estimate}), and the error related to the approximation in the stochastic domain. The latter is an interpolation error, which is the main focus of discussion in this section.
	In this context, we demonstrate the quasi-optimality of the interpolant approximation in the polynomial space $\mathcal{P}_{\boldsymbol{p}}(\boldsymbol{\Gamma})$ in Lemma~\ref{lem:optimality-interpolation}. 
	To achieve this, we impose certain restrictions on the growth of $\rho$ at infinity, as presented in Assumption~\ref{assum:convergence-A-f}~\textit{\ref{assum:upper-bound-density}}. Specifically, we assume that $\rho$ behaves like a Gaussian probability density at infinity. Additionally, to apply the upper bounds provided by Lemmas~\ref{lem:bonded-domain-min-approx} and \ref{lem:unbonded-domain-min-approx} for the interpolant approximation, it is necessary to ensure two conditions:
	\begin{enumerate}
		\item 
		For unbounded $\boldsymbol{\Gamma}$, the growth of the solution $u$ at infinity is at most exponential.
		\item 
		The solution $u$ admits an analytic extension with respect to each variable $y_n$, for $n=1, \cdots, N$.
	\end{enumerate}
	To control the growth of $u$ at infinity by regulating the growth of $f$ , we rely on Assumption~\ref{assum:convergence-A-f}~\textit{\ref{assum:f-stochastic-domin}}, and to guarantee the analytic extension of $u$, we rely on Assumption~\ref{assum:convergence-A-f}~\textit{\ref{assum:boundedness-derivative-random-parameter}}.
	For the remainder of this work, we assume that $V$ is a Hilbert functional space.
	
	\vspace*{1mm}
	\begin{definition}[{\cite[(3.1)]{Babuska-Nobile-Tempone}}]
		We introduce a weight $\sigma(\boldsymbol{y})= \prod _{n=1}^N \sigma_n(y_n) \leq 1$ where
		\begin{equation}
			\label{ieq:sigma-weight}
			\sigma_n(y_n):=
			\begin{cases}
				1 & \text{if $\Gamma_n$ is bounded},
				\\
				e^{-\alpha_n \lvert y_n \rvert }
				& \text{if $\Gamma_n$ is unbounded},
			\end{cases}
		\end{equation}
		for some $\alpha_n>0$ and then define the functional space
		\begin{equation}
			\label{def-C-sigmma-space}
			C_\sigma^0(\boldsymbol{\Gamma}; V):=
			\left\lbrace 
			\hspace*{-0.5mm}
			v:\boldsymbol{\Gamma} \rightarrow V
			%, ~
			{\big \vert}~ 
			v~ \text{is continuous in } \boldsymbol{y}, 
			\hspace*{0.5mm}
			\left\| v \right\|_{C_\sigma^0(\boldsymbol{\Gamma}; V)}
			\hspace*{-0.5mm}
			:=
			%{\color{red}
			\max_{\boldsymbol{y} \in \boldsymbol{\Gamma}} 
			 \sigma(\boldsymbol{y}) \lVert v(\boldsymbol{y}) \rVert_{V} 
			%}  
			 < \infty
			\hspace*{-0.5mm}
			\right\rbrace.
		\end{equation}
	\end{definition}
	
	\vspace*{1mm}
	%\begin{definition}
		For any $\boldsymbol{M}:\boldsymbol{\Gamma} \times \mathcal{D} \rightarrow \mathbb{R}^{d \times d}$ with
		\begin{equation}
			\boldsymbol{M}(\boldsymbol{y}, \boldsymbol{x}) = 
			\begin{bmatrix}
				m_{1,1}(\boldsymbol{y}, \boldsymbol{x}) & \cdots & m_{1,d}(\boldsymbol{y}, \boldsymbol{x}) \\
				\vdots & \ddots & \vdots \\
				m_{d,1}(\boldsymbol{y}, \boldsymbol{x}) & \cdots & m_{d,d}(\boldsymbol{y}, \boldsymbol{x})
			\end{bmatrix},
		\end{equation}
		we define the $k$-order partial derivative of 
		$\boldsymbol{M}$ with respect to the $y_n$, for each $k \in \mathbb{N}$ as 
		\begin{equation}
			\label{def-partial-derveitive-matrix}
			\partial^k_{y_n} \boldsymbol{M}:= 
			\begin{bmatrix}
				\partial_{y_n}^k(m_{1,1}) & \cdots &
				\partial_{y_n}^k(m_{1,d})
				\\
				\vdots & \ddots & \vdots
				\\
				\partial_{y_n}^k(m_{d,1}) & \cdots & \partial_{y_n}^k(m_{d,d}) 
			\end{bmatrix}.
		\end{equation}
	%\end{definition}
	
	In the following, we denote $\bullet_n$ as a quantity related to direction $y_n$ and $\bullet_n^*$ as the analogous quantity related to all other directions $y_j$ where $j\neq n$.
	
	\vspace*{1mm}
	\begin{assumption}[{ \cite[Assumption~2 and Lemma~3.2]{Babuska-Nobile-Tempone}}]
		\label{assum:convergence-A-f}
		The joint probability density $\rho$ and the random data $f$ and $\boldsymbol{A}$ satisfy the following conditions:
		\begin{enumerate}[label=\textnormal{(\alph*)}]
			\item
			\label{assum:upper-bound-density}
			There exists a constant $C_{\ref{eq:density-func-assum}}>0$ such that 
			\begin{equation}
				\label{eq:density-func-assum}
				\rho(\boldsymbol{y}) 
				\leq
				C_{\ref{eq:density-func-assum}}  
				e^{- \sum_{n=1}^N (\delta_n y_n)^2 },
				~~ \forall \boldsymbol{y} \in \boldsymbol{\Gamma},
			\end{equation}
			where $\delta_n>0$ if $\Gamma_n$ is unbounded, and $\delta_n=0$ otherwise.
			\item 
			\label{assum:f-stochastic-domin}
			$f \in C^0_{\sigma}\left( \boldsymbol{\Gamma}; L^2(\mathcal{D}) \right)$.
			\item 
			\label{assum:boundedness-derivative-random-parameter}
			There exist constants $\eta_n < \infty$ such that for every $\boldsymbol{y} = (y_n, \boldsymbol{y}^*_n) \in \boldsymbol{\Gamma}$ and $k \in \mathbb{N}$,
			\begin{equation}
				\label{assum:convergence}
				\frac{\left\|  \partial^k_{y_n} f(\boldsymbol{y}) \right\| _{L^2(\mathcal{D})}}{1+ \left\|  f(\boldsymbol{y}) \right\| _{L^2(\mathcal{D})}}
				\leq
				\eta_n^k k!,
				\quad
				\text{and}
				\quad
				\left\| 
				\left( \partial^k_{y_n} \boldsymbol{A} (\boldsymbol{y}) \right) 
				 \boldsymbol{A}^{-1}(\boldsymbol{y})
				\right\| _{L^\infty(\mathcal{D})}
				\leq
				\eta_n^k k!
				.
			\end{equation}
		\end{enumerate}
	\end{assumption}
	
	\vspace*{1mm}
	\begin{remark}[setting auxiliary density function {\cite[(3.3)]{Babuska-Nobile-Tempone}}]
		Based on Assumption~\ref{assum:convergence-A-f}~\textit{\ref{assum:upper-bound-density}}, we choose a suitable auxiliary density function $\hat{\rho}(\boldsymbol{y})= \prod_{n=1}^N \hat{\rho}(y_n)$ such that for each $n= 1, \dots,N$ there exist $C^{m}_n, C^{M}_n>0$ satisfying
		\begin{equation}
			\label{ineq:auxiliary-density}
			C^{m}_n e^{-(\delta_ny_n)^2}
			\leq
			\hat{\rho}_n(y_n)
			\leq
			C^{M}_n e^{-(\delta_ny_n)^2}, 
			\quad
			\forall
			y_n \in \Gamma_n.
		\end{equation}
		This choice satisfies the requirement~(\ref{def:join-rho}); indeed, 
		\begin{equation}
			\label{ineq:upper-bound-rho-on-hat}	
		\left\| \dfrac{\rho}{\hat{\rho}}\right\| _{L^\infty(\boldsymbol{\Gamma})} \leq \dfrac{C_{\ref{eq:density-func-assum}}}{\prod_{n=1}^{N} C^{m}_n}.
		\end{equation}
	\end{remark}
	
	\vspace*{1mm}
	\begin{remark}[regularity of strong solution {\cite[Lemma~3.1]{Babuska-Nobile-Tempone}}]
		\label{rem:regularity-w.r-paramter}
		Assumption~\ref{assum:convergence-A-f}~\textit{\ref{assum:f-stochastic-domin}}, together with $\boldsymbol{A} \in C^0_{\text{loc}}\left( \boldsymbol{\Gamma} ; L^{\infty}(\mathcal{D}, \mathbb{R}^{d \times d}) \right) $, immediately implies that the strong solution of (\ref{eq:nondivergence-random-variable}) satisfies $u \in C_\sigma^0( \boldsymbol{\Gamma}; H^2) $.
	\end{remark}
	
	\vspace*{1mm}
	\begin{definition}
		\label{def:space-versigma}
		We define the functional space $C_{\varsigma}^0\left( \boldsymbol{\Gamma}, V \right) $ as same as $C_\sigma^0(\boldsymbol{\Gamma}; V)$, except that 
		the wight $\sigma(\boldsymbol{y}) = \prod _{n=1}^N \sigma_n(y_n)$ is replaces by $\varsigma(\boldsymbol{y}) = \prod_{n=1}^{N}\varsigma_n(y_n)$, where 
		\begin{equation}
			\label{eq:varsigma-weight}
			\varsigma_n(y_n):=
			\begin{cases}
				1 & \text{if $\Gamma_n$ is bounded},
				\\
				e^{-\frac{(\delta_ny_n)^2}{4}} 
				& \text{if $\Gamma_n$ is unbounded},
			\end{cases}	
		\end{equation}
		and $\delta_n$ are specified from (\ref{eq:density-func-assum}).
	\end{definition}
	
	\vspace*{1mm}
	The motivation for defining the additional weight (\ref{eq:varsigma-weight}) and, consequently, the space $C_{\varsigma}^0\left( \boldsymbol{\Gamma}, V \right)$ is that the upper bound in Lemma~\ref{lem:unbonded-domain-min-approx} is expressed in terms of this weight.
	%	{\color{red}To establish a priori estimates for the total error
		%	\begin{equation*}
			%	\lVert (u, \boldsymbol{g}) - (u_{p,\mathbb{U}}, \boldsymbol{g}_{p, \mathbb{G}})\rVert_{\left( L^2_\rho(\boldsymbol{\Gamma}) \otimes H^1(D)\right)  \times \left( L^2_\rho(\boldsymbol{\Gamma}) \otimes H^1(D; \mathbb{R}^d) \right) },
			%\end{equation*}
			%	we must first derive and prove several supporting lemmas. These lemmas will serve as essential tools for the error analysis and ensure the validity of the estimates.}
			
			\vspace*{1mm}
		\begin{lemma}[continuously embedding]
			\label{lem:countiuity-embedding}
			Under Assumption~\ref{assum:convergence-A-f}~\textit{\ref{assum:upper-bound-density}}, the following continuous embedding holds
			\begin{equation}
				\label{ineq:countiuity-embedding}
				C_\sigma^0(\boldsymbol{\Gamma}; V) 
				\subseteq
				C_\varsigma^0(\boldsymbol{\Gamma}; V) 
				\subseteq
				L^2_{\hat{\rho}}(\boldsymbol{\Gamma}; V)
				\subseteq
				L^2_{\rho}(\boldsymbol{\Gamma}; V).
			\end{equation}
		\end{lemma}
		
		\vspace*{1mm}
		This inclusion is a more complete form of what is presented in {\cite[Section 3]{Babuska-Nobile-Tempone}}, as it incorporates the intermediate space $C_\varsigma^0(\boldsymbol{\Gamma}; V) $, which was not included there.
		
		\vspace*{1mm}
		\begin{proof}
			We will demonstrate the inclusions in (\ref{ineq:countiuity-embedding}) step by step, starting from the last inclusion $	L^2_{\hat{\rho}}(\boldsymbol{\Gamma}; V) \subseteq L^2_{\rho}(\boldsymbol{\Gamma}; V)$ and proceeding to the first inclusion $C_\sigma^0(\boldsymbol{\Gamma}; V) \subseteq C_\varsigma^0(\boldsymbol{\Gamma}; V) $.
			%\\
			From (\ref{ineq:upper-bound-rho-on-hat}), we have
			\begin{equation*}
				\left\| v\right\|_{L^2_{\rho}(\boldsymbol{\Gamma}; V)}^2
				\leq
				\left\| \dfrac{\rho}{\hat{\rho}}\right\| _{L^\infty(\boldsymbol{\Gamma})} 
				\left\| v\right\|_{L^2_{\hat{\rho}}(\boldsymbol{\Gamma}; V)}^2
				\leq
				\dfrac{C_{\ref{eq:density-func-assum}}}{\prod_{n=1}^N C^{m}_n}
				\left\| v\right\|_{L^2_{\hat{\rho}}(\boldsymbol{\Gamma}; V)}^2.
			\end{equation*}
			We now, aim to demonstrate the inclusion $C_\varsigma^0(\boldsymbol{\Gamma}; V) \subseteq L^2_{\hat{\rho}}(\boldsymbol{\Gamma}; V)$.
			\begin{equation*}
				\left\| v\right\|_{L^2_{\hat{\rho}}(\boldsymbol{\Gamma}; V)}^2 = \int_{\boldsymbol{\Gamma}} \left\|  v(\boldsymbol{y}) \right\|_V^2  \hat{\rho}(\boldsymbol{y}) d\boldsymbol{y} 
				\leq
				\left\| v\right\|^2_{C_\varsigma^0(\boldsymbol{\Gamma}; V)} \int_{\boldsymbol{\Gamma}} \dfrac{\hat{\rho}(\boldsymbol{y})}{\varsigma^2(\boldsymbol{y})} d\boldsymbol{y} 
				\leq
				\prod_{n=1}^N C_n 
				\left\|v\right\|^2_{C_\varsigma^0(\boldsymbol{\Gamma}; V)},
			\end{equation*}
			%	where for bounded $\Gamma_n$, $M_n := C_{\text{max}}^n \left| \Gamma_n\right| $ and for unbounded $\Gamma_n$, $M_n := C_{\text{max}}^n e^{(\alpha_n/\delta_n)^2} \dfrac{\sqrt{\pi}}{\delta_n}$, because 
			%	\begin{multline*}
				%		\int_{\Gamma_n}\dfrac{\hat{\rho}_n(y_n)}{\sigma_n^2(y_n)}~ dy_n = \int_{\Gamma_n} \dfrac{\hat{\rho}_n(y_n)}{e^{-2\alpha_n\left| y_n\right| }} ~ dy_n
				%		\leq
				%		C_{\text{max}}^n \int_{\Gamma_n} e^{-(\delta_n y_n)^2 + 2\alpha_n |y_n|}~ dy_n
				%		\\
				%		=
				%		C_{\text{max}}^n e^{{\left( \alpha_n / \delta_n\right) }^2} \int_{\Gamma_n} e^{-\left( \delta_n 	|y_n| - \left( \alpha_n / \delta_n\right) \right)^2 } ~dy_n
				%		\leq
				%		C_{\text{max}}^n e^{{\left( \alpha_n / \delta_n\right) }^2} \dfrac{\sqrt{\pi}}{\delta_n}.
				%	\end{multline*}
			where for bounded $\Gamma_n$, $C_n := C^{M}_n \left| \Gamma_n\right| $ and for unbounded $\Gamma_n$, $C_n := C^{M}_n \frac{\sqrt{2\pi}}{\delta_n}$, because 
			\begin{equation*}
				\int_{\Gamma_n}\dfrac{\hat{\rho}_n(y_n)}{\varsigma_n^2(y_n)} dy_n 
				= 
				\int_{\Gamma_n} \dfrac{\hat{\rho}_n(y_n)}{e^{-(\delta_n y_n)^2/2 }}  dy_n
				\leq
				C^{M}_n \int_{\Gamma_n} e^{-\frac{(\delta_n y_n)^2}{2}} dy_n
				\\
				\leq
				C^{M}_n \dfrac{\sqrt{2\pi}}{\delta_n}.
			\end{equation*}
			Finally, we show that $C_\sigma^0(\boldsymbol{\Gamma}; V) \subseteq C_\varsigma^0(\boldsymbol{\Gamma}; V) $.
			\begin{equation*}
				\left\| v  \right\|_{C_{\varsigma}^0(\boldsymbol{\Gamma}; V)}
				=
				%{\color{red}
				\max_{\boldsymbol{y} \in \boldsymbol{\Gamma}}  \varsigma(\boldsymbol{y}) \left\| v(\boldsymbol{y}) \right\|_V  
				%}
				\leq
				%{\color{red}
				\max_{\boldsymbol{y} \in \boldsymbol{\Gamma}}   \dfrac{\varsigma(\boldsymbol{y})}{\sigma(\boldsymbol{y})} 
				\left\| v  \right\|_{C_{\sigma}^0(\boldsymbol{\Gamma}; V)}
				%}
				\leq
				\prod_{n=1}^{N}\tilde{C}_n 
				\left\| v  \right\|_{C_{\sigma}^0(\boldsymbol{\Gamma}; V)},
			\end{equation*}
			where for bounded $\Gamma_n$, since $\sigma_n(y_n)=\varsigma_n(y_n)=1 $, thus $\tilde{C}_n=1$ and for unbounded $\Gamma_n$, $\tilde{C}_n= e^{\left( \frac{\alpha_n}{\delta_n} \right)^2 }$ because
			\begin{equation*}
				%{\color{red}
				\max_{y_n \in \Gamma_n}  e^{-\frac{(\delta_n y_n)^2}{4}  + \alpha_n \left| y_n\right| } 
				=
				\max_{y_n \in \Gamma_n}   e^{-\left( \frac{\delta_n \left| y_n\right| }{2}  - \frac{\alpha_n}{\delta_n} \right)^2 + \left( \frac{\alpha_n}{\delta_n}\right)^2  }  
				%}
				\leq
				e^{\left( \frac{\alpha_n}{\delta_n} \right)^2 }.
			\end{equation*}
		\end{proof}
		%	In lemmas~\ref{lem:bonded-domain-min-approx} and \ref{lem:unbonded-domain-min-approx}, which provide upper bounds for interpolation error in bounded and unbounded domains, respectively, these lemmas establish such upper bounds for functions that admit analytic extensions. The upper bounds are determined by the diameter of the complex region where the function is analytic. To proceed, in the next lemma, we verify the analytic extension of the solution $u$ for each variable $y_n$, where $n=1, \cdots,N$.
		
		\vspace*{1mm}
		\begin{lemma}[analytic extension]
			\label{lem:analytic-exte}
			Under Assumption~\ref{assum:convergence-A-f}~\textit{\ref{assum:boundedness-derivative-random-parameter}}, the solution $u(y_n, \boldsymbol{y}_n^*, \boldsymbol{x})$ as a function of $y_n$, $u:\Gamma_n \rightarrow C_{\sigma_n^*}\left( \boldsymbol{\Gamma}_n^*; H^2 \right)$ admits an analytic extension $u(z, \boldsymbol{y}_n^*, \boldsymbol{x})$, $z \in \mathbb{C}$, in the region of the complex plane 
			\begin{equation*}
				\Sigma(\Gamma_n; \tau_n):=\left\lbrace 
				z \in \mathbb{C},~ \operatorname{dist}(z, \Gamma_n) \leq \tau_n
				\right\rbrace,  \quad 
				0 < \tau_n < \frac{1}{2\eta_n}.
			\end{equation*}
			Moreover, $\forall z \in \Sigma(\Gamma_n; \tau_n)$,
			\begin{equation}
				\label{analytic-solution-upper-bound}
				\left\| 
				\sigma_n(Re~z) u(z)
				\right\|_{C^0_{\sigma_n^*}\left( \boldsymbol{\Gamma}_n^*; H^2(\mathcal{D}) \right) } 
				\leq
				%\\
				\dfrac{e^{\alpha_n \tau_n}}{\sqrt{C_{\ref{ineq:coercivity-a}}}(1-2\tau_n \eta_n)} \left( 
				2 \left\| f \right\|_{C^0_{\sigma}\left( \boldsymbol{\Gamma}; L^2(\mathcal{D}) \right) } +1
				\right).
			\end{equation}
		\end{lemma}
		
		\vspace*{1mm}
		This lemma is an adaptation of {\cite[Lemma~3.2]{Babuska-Nobile-Tempone}} to equation (\ref{eq:nondivergence-random-variable}).
		
		\vspace*{1mm}
		\begin{proof}
			For any $k \in \mathbb{N}$, taking the $k$-the derivative with respect to $y_n$ of both sides of equation (\ref{eq:nondivergence-random-variable}) results in the following equality
			\begin{equation*}
				\boldsymbol{A}(\boldsymbol{y}): \mathrm{D}^2\left( \partial_{y_n}^k u(\boldsymbol{y}) \right) 
				= 
				- \sum_{j = 1}^{k} \binom{k}{j} \partial_{y_n}^j \boldsymbol{A}(\boldsymbol{y}) : \mathrm{D}^2\left( \partial_{y_n}^{k-j} u(\boldsymbol{y}) \right) 
				+
				\partial_{y_n}^kf(\boldsymbol{y}).
			\end{equation*}
			Hence, this implies that
			\begin{multline*}
				\left\| 
				\boldsymbol{A}(\boldsymbol{y}): \mathrm{D}^2\left( \partial_{y_n}^k u(\boldsymbol{y}) \right)
				\right\|_{L^2(\mathcal{D})} 
				\\
				\leq
				%\\
				\sum_{j = 1}^{k} \binom{k}{j}
				\left\| 
				\left( \partial_{y_n}^j \boldsymbol{A}(\boldsymbol{y}) \right)  \boldsymbol{A}^{-1}(\boldsymbol{y}) 
				\right\|_{L^\infty(\mathcal{D})} 
				\left\| 
				\boldsymbol{A}(\boldsymbol{y}): \mathrm{D}^2\left( \partial_{y_n}^{k-j} u(\boldsymbol{y}) \right) 
				\right\|_{L^2(\mathcal{D})} 
				+
				\left\| 
				\partial_{y_n}^kf(\boldsymbol{y})
				\right\|_{L^2(\mathcal{D})}.  
			\end{multline*}
			Setting $R_k(\boldsymbol{y}) = {\left\| \boldsymbol{A}(\boldsymbol{y}): \mathrm{D}^2\left( \partial_{y_n}^k u(\boldsymbol{y}) \right) \right\|_{L^2(\mathcal{D})}}/{k!} $ and using the upper bounds (\ref{assum:convergence}), we get the recursive inequality 
			\begin{equation*}
				R_k(\boldsymbol{y}) \leq \sum_{j=1}^{k} \eta_n^j R_{k-j}(\boldsymbol{y}) + \eta_n^k\left( 1+ \left\| f(\boldsymbol{y}) \right\|_{L^2(\mathcal{D})}  \right). 
			\end{equation*}
			This leads to 
			\begin{equation*}
				R_k (\boldsymbol{y}) \leq \left( 2\eta_n \right) ^k R_0 (\boldsymbol{y}) + \left( 1+ \left\| f(\boldsymbol{y}) \right\|_{L^2(\mathcal{D})}  \right) \eta_n^k \sum_{j=1}^{k}2^{j-1}. 
			\end{equation*} 
			Since $R_0(\boldsymbol{y}) =\left\| \boldsymbol{A}(\boldsymbol{y}) : \mathrm{D}^2u(\boldsymbol{y}) \right\|_{L^2(\mathcal{D})} = \left\| f(\boldsymbol{y}) \right\|_{L^2(\mathcal{D})}$, it follows that
			\begin{equation*}
				\dfrac{\left\| \boldsymbol{A}(\boldsymbol{y}): \mathrm{D}^2\left( \partial_{y_n}^k u(\boldsymbol{y}) \right) \right\|_{L^2(\mathcal{D})}}{k!}
				\leq
				(2\eta_n)^k \left( 2\left\| f(\boldsymbol{y}) \right\|_{L^2(\mathcal{D})} +  1 \right) .
			\end{equation*}
			Thanks to Lemma~\ref{lem:well-posedness-a-y}, there exists some constant $C_{\ref{ineq:coercivity-a}}>0$ such for every $v \in H^2$, $C_{\ref{ineq:coercivity-a}} \left\| v\right\|_{H^2(\mathcal{D})} \leq  \left\|\boldsymbol{A}(\boldsymbol{y}): \mathrm{D}^2v\right\|_{L^2(\mathcal{D})} $,
			\begin{equation*}
				\dfrac{\left\|  \partial_{y_n}^k u(\boldsymbol{y}) \right\|_{H^2(\mathcal{D})}}{k!}
				\leq
				(2\eta_n)^k \frac{ \left( 2\left\| f(\boldsymbol{y}) \right\|_{L^2(\mathcal{D})} +  1 \right)}{\sqrt{C_{\ref{ineq:coercivity-a}}}} .
			\end{equation*}
			We now define, for every $y_n \in \Gamma_n$ the power series $u: \mathbb{C} \rightarrow C^0_{\sigma_n^*}\left( \boldsymbol{\Gamma}_n^*, H^2(\mathcal{D}) \right) $ as 
			\begin{equation*}
				u\left( z, \boldsymbol{y}_n^*, \boldsymbol{x} \right) = \sum_{k = 0}^{\infty} \dfrac{(z - y_n)^k}{k!} \partial_{y_n}^k \left( y_n, \boldsymbol{y}_n^*, \boldsymbol{x} \right) .
			\end{equation*}
			Therefore,
			\begin{multline*}
				\sigma_n(y_n) \left\| u(z) \right\| _{C^0_{\sigma_n^*}\left( \boldsymbol{\Gamma}_n^*; H^2(\mathcal{D}) \right)}  
				\leq
				\sum_{k=0}^{\infty} \dfrac{\left| z - y_n \right|^k }{k!}\sigma_n(y_n) \left\| \partial^k_{y_n} u(y_n) \right\|_{C^0_{\sigma_n^*}\left( \boldsymbol{\Gamma}_n^*; H^2(\mathcal{D}) \right)}  
				\\ 
				\leq
				(C_{\ref{ineq:coercivity-a}})^{-1/2} 
				%{\color{red}
				\max_{y \in \Gamma_n} \sigma_n(y) \left( 2 \left\| f(y)\right\| _{C^0_{\sigma_n^*}(\boldsymbol{\Gamma}_n^*; L^2(\mathcal{D}))} +1 \right)
				%} 
				\sum_{k=0}^{\infty} \left( 2\eta_n \left| z-y_n \right|  \right)^k
				\\
				\leq
				\left( C_{\ref{ineq:coercivity-a}} \right)^{-1/2}\left( 2 \left\| f\right\|_{C^0_\sigma(\boldsymbol{\Gamma}; L^2(\mathcal{D}))} +1  \right)
				\sum_{k=0}^{\infty} \left(2\eta_n \left| z - y_n \right|  \right)^k .   
			\end{multline*}
			The series converges for all $z \in \mathbb{C}$ such that $\left|  z- y_n \right| \leq \tau_n <{1}/{(2 \eta_n)}$. Within the ball $\left| z - y_n \right| \leq \tau_n$, we have $\sigma_n(Re~ z) \leq e^{\alpha_n \tau_n}\sigma_n(y_n)$. Therefore, we have
			\begin{equation*}
				\sigma_n(Re~ z) \left\| u(z) \right\|_{C^0_{\sigma_n^*}\left( \boldsymbol{\Gamma}_n^*; H^2(\mathcal{D}) \right)}  
				\leq \dfrac{e^{\alpha_n \tau_n}}{\sqrt{C_{\ref{ineq:coercivity-a}}}(1-2\tau_n \eta_n)}
				\left( 2 \left\| f\right\|_{C^0_\sigma(\boldsymbol{\Gamma}; L^2(\mathcal{D}))} + 1 \right). 
			\end{equation*}
			As the power series converges for every $y_n \in \Gamma_n$, by a continuation argument, the function $u$ can be extended analytically on the whole region $\Sigma(\Gamma_n; \tau_n)$ and estimate~(\ref{analytic-solution-upper-bound}) follows.
		\end{proof}
		
		\vspace*{1mm}
		The $N$-dimensional interpolation operator can be viewed as a composition of $N$ one-dimensional interpolants. 
		Consequently, we proceed by discussing the error bound for the one-dimensional interpolation operator. 
		In Lemmas~\ref{lem:continuity-interpolation-op}–\ref{lem:unbonded-domain-min-approx}, we recall some established results from approximation theory for functions defined on a one-dimensional domain with values in a Hilbert space $V$. 
		To this end, we use the notations introduced earlier for the one-dimensional domain $\Gamma$. 
		Let $\rho: \Gamma \rightarrow \mathbb{R}^+ \cup \left\lbrace 0\right\rbrace $ be a weight satisfying $\rho(y) \leq C^{M}e^{-(\delta y)^2}$ for all $ y \in \Gamma$, where $C^{M} > 0$ and $\delta>0$ if $\Gamma$ is unbounded, and $\delta = 0$ otherwise. 
		Corresponding to this $\delta$, we define a weight function $\varsigma(y)= e^{-{(\delta y)^2}/{4}}$. 
		Additionally, we introduce another weight function $\sigma$ defined as $\sigma(y)=e^{-\alpha |y|}$, where $\alpha>0$ if $\Gamma$ is unbounded, and $\alpha = 0$ otherwise. The space $C^0_\sigma\left( \Gamma; V \right)$ is then defined similarly to (\ref{def-C-sigmma-space}), but in the one-dimensional case where $N=1$. 
		The space $C^0_\varsigma\left( \Gamma; V \right)$ is defined analogously, with $\sigma$ replaced by $\varsigma$. We let $y_j \in \Gamma$ for $j= 1,\dots, p+1$, denote the zeros of a polynomial of degree $p+1$ orthogonal to the space $\mathcal{P}_p(\Gamma)$ with respect to the weight $\rho$. 
		For a continuous function $v:\Gamma \subset \mathbb{R} \rightarrow V $, the Lagrange interpolation $\mathcal{I}_p$ is defined as $\mathcal{I}_pv(y):=\sum_{j=1}^{p+1}v(y_j) \ell_j (y)$, where $\ell_j (y)$ are the Lagrange basis polynomials. 
		The weights of the Gaussian quadrature formula are defined as  $w_j:= \int_{\Gamma}\ell_j(y) \rho(y)~dy = \int_{\Gamma}\ell_j^2(y) \rho(y)~dy$. 
		Using these notations and arguments similar to the proof of Lemma~\ref{lem:countiuity-embedding}, in the case $N=1$ the following inclusions hold
		\begin{equation}
			\label{subset:inclusions}
			C^0_\sigma\left( \Gamma; V \right)
			\subseteq
			C^0_\varsigma\left( \Gamma; V \right)
			\subseteq
			L^2_\rho \left( \Gamma; V \right) .
		\end{equation}	

			\vspace*{1mm}
			\begin{lemma}[continuity of interpolation]
				\label{lem:continuity-interpolation-op}
				The operator $\mathcal{I}_p:C^0_\varsigma\left( \Gamma; V \right) \rightarrow L^2_\rho\left( \Gamma; V \right)$ is continuous; i.e., the exists a constant $C_{\ref{ineq:continuity-interpolation-op}}>0$ such that for any $v \in C^0_\varsigma\left( \Gamma; V \right)$ the following holds
				\begin{equation}
					\label{ineq:continuity-interpolation-op}
					\left\| \mathcal{I}_p v\right\|_{L_\rho^2(\Gamma; V)} 
					\leq
					C_{\ref{ineq:continuity-interpolation-op}}
					\left\| v \right\|_{C_\varsigma^0(\Gamma; V)}.
				\end{equation}
			\end{lemma}
			
			\vspace*{1mm}
			This lemma is a variation of Lemma 4.2 from \cite{Babuska-Nobile-Tempone}, where the functional space $V$ is considered a Hilbert space instead of a Banach space, and the weight function is $\varsigma$ instead of $\sigma$.
			
			\vspace*{1mm}
			\begin{proof}
				For any $v \in C^0_\varsigma\left( \Gamma; V \right)$, where $V$ is a Hilbert space equipped with the inner product $(\cdot, \cdot)_V$, we have
				\begin{align*}
					\left\| \mathcal{I}_p v\right\|_{L_\rho^2(\Gamma; V)}^2 
					&= 
					\bigintss_{\Gamma} \left\| \sum_{j=1}^{p+1} v(y_j) \ell_j(y) \right\|_V^2 \rho(y)dy 
					\\ 
					&=
					\bigintsss_{\Gamma} \left( \sum_{j=1}^{p+1} v(y_j) \ell_j(y), \sum_{j=1}^{p+1} v(y_j) \ell_j(y) \right)_V \rho(y)dy. 
				\end{align*}
				Thanks to the orthogonality property 
				$\int_{\Gamma}\ell_i(y)\ell_j(y)\rho(y)~dy = 0,~ \forall i\neq j$, we have
				\begin{align*}
					\left\| \mathcal{I}_p v\right\|_{L_\rho^2(\Gamma; V)}^2 
					&= 
					\bigintssss_{\Gamma} \sum_{j=1}^{p+1} \left\| v(y_j)\right\|_V^2  \ell_j^2(y) \rho(y)dy
					\\
					&\leq
					%{\color{red}
					\max_{y \in \left\lbrace y_1, \dots, y_{p+1} \right\rbrace} 
					\left\| v(y)\right\|_V^2 \varsigma^2(y)
					%}
					\sum_{j=1}^{p+1} \bigintssss_{\Gamma} \frac{\ell_j^2(y) \rho(y)}{\varsigma^2(y_j)}dy.
				\end{align*}
				Then, Remark~\ref{rem:equivalence-weights} implies that
				\begin{equation*}
					\left\| \mathcal{I}_p v\right\|_{L_\rho^2(\Gamma; V)}^2 
					\leq
					\left\| v \right\|_{C_\varsigma^0(\Gamma; V)}
					\sum_{j=1}^{p+1} \dfrac{w_j}{\varsigma^2(y_j)}. 
				\end{equation*}
				For the case of bounded $\Gamma$, since $\varsigma(y)=1$ and $\sum_{k=1}^{p+1} w_k = 1$, the result follows. For unbounded $\Gamma$, sine $\rho(y) \leq  C^M e^{-(\delta y)^2}$ and all the even moments $c_{2m} = \int_{\Gamma}y^{2m} \rho(y)dy$ are uniformly bounded, using a result from \cite{Uspensky} leads to
				\begin{equation*}
					\sum_{j=1}^{p+1} \dfrac{w_j}{\varsigma^2(y_j)}
					\xrightarrow[]{p\rightarrow \infty}
					\int_{\Gamma}\frac{\rho(y)}{\varsigma^2(y)}
					dy \leq 
					C^M \int_{\Gamma} e^{-(\delta y)^2/2}dy
					\leq 
					C^M \frac{\sqrt{2 \pi}}{\delta},
				\end{equation*} 
				and we can conclude that (\ref{ineq:continuity-interpolation-op}) holds.
			\end{proof}
			
			\vspace*{1mm}
			\begin{lemma}[quasi-optimality {\cite[Lemma~4.3]{Babuska-Nobile-Tempone}}]
				\label{lem:optimality-interpolation}
				There exists a constant $C_{\ref{ineq:optimality-interpolation}}>0$ such that for every $v \in C^0_\sigma\left( \Gamma; V \right)$, the following holds
				\begin{equation}
					\label{ineq:optimality-interpolation}
					\left\| v - \mathcal{I}_p v \right\|_{L_\rho^2(\Gamma; V)}
					\leq
					C_{\ref{ineq:optimality-interpolation}}
					%{\color{red}
					\inf_{w \in \mathcal{P}_p(\Gamma) \otimes V}
					\left\| v - w \right\|_{C^0_\varsigma\left( \Gamma; V \right)}
					%}
					.  
				\end{equation}
			\end{lemma}
			
			\vspace*{1mm}
			\begin{proof}
				For every $w \in \mathcal{P}_p(\Gamma) \otimes V$, since $\mathcal{I}_p w = w$, from (\ref{subset:inclusions}) and Lemma~\ref{lem:continuity-interpolation-op} we have
				\begin{equation*}
					\left\| v - \mathcal{I}_p v \right\|_{L^2_\rho(\Gamma; V)} 
					\leq 
					\left\| v - w \right\|_{L^2_\rho(\Gamma; V)} 
					+
					\left\| \mathcal{I}_p(w - v) \right\|_{L^2_\rho(\Gamma; V)}
					\leq
					C_{\ref{ineq:optimality-interpolation}} \left\|  v - w  \right\|_{C^0_\varsigma\left( \Gamma; V \right)}.
				\end{equation*}
			\end{proof}
			
			\vspace*{1mm}
			\begin{lemma}[interpolation error bound in a bounded domain {\cite[Lemma~4.4]{Babuska-Nobile-Tempone}}]
				\label{lem:bonded-domain-min-approx}
				A function $v \in C^0(\Gamma; V)$ which admits an analytic extension in the complex region $\Sigma \left( \Gamma; \tau \right) = \left\lbrace z \in \mathbb{C},~\operatorname{dist}(z, \Gamma) \leq \tau \right\rbrace $ for some $\tau>0$, satisfies
				\begin{equation}
					\label{ineq:bonded-domain-min-approx}
					%{\color{red}
					\min_{\phi \in \mathcal{P}_p(\Gamma) \otimes V}
					\left\| v - \phi \right\|_{C^0\left( \Gamma; V \right)}
					%}
					\leq
					\dfrac{2}{\varrho - 1}e^{-p \log(\varrho)} 
					%{\color{red} 
					\max_{z \in \Sigma \left( \Gamma; \tau \right)} \left\| v(z)\right\|_V
					%}
					, 		
				\end{equation}
				where $1<\varrho := {2\tau}/{\left| \Gamma \right| } + \sqrt{1 + {4 \tau^2}/{\left|  \Gamma\right| ^2}}$.
			\end{lemma}
			
			\vspace*{1mm}
			\begin{proof}
				Refer to Lemma~4.4 of \cite{Babuska-Nobile-Tempone}.
			\end{proof}
			
			\vspace*{1mm}
			\begin{lemma}[interpolation error bound in an unbounded domain {\cite[Lemma~4.6]{Babuska-Nobile-Tempone}}]
				\label{lem:unbonded-domain-min-approx}
				Let $v \in C^0_\sigma (\mathbb{R}; V)$ be a function which admits an analytic extension in the complex strip  $\Sigma \left( \Gamma; \tau \right) = \left\lbrace z \in \mathbb{C},~ \operatorname{dist}(z, \Gamma) \leq \tau \right\rbrace $ for some $\tau>0$. Moreover, for each $z = (y + \boldsymbol{i} w) \in \Sigma(\mathbb{R}; \tau)$, it satisfies
				\begin{equation*}
					\sigma(y) \left\| v(z)\right\|_{V} \leq C_v(\tau),
				\end{equation*}
				where $\boldsymbol{i}$ denotes the imaginary unit.
				Then, for any $\delta>0$, there exists a constant $C_{\ref{ineq:unbonded-domain-min-approx}}>0$, independent of $p$, such that
				\begin{equation}
					\label{ineq:unbonded-domain-min-approx}
					%{\color{red}
					\min_{\phi \in \mathcal{P}_p \otimes V} \max_{y \in \mathbb{R}} \left\| v(y) - \phi(y) \right\|_V  e^{-\frac{(\delta y)^2}{4}}
					%}
					\leq
					C_{\ref{ineq:unbonded-domain-min-approx}}
					O(\sqrt{p})  e^{-\tau \delta\sqrt{p}}.
				\end{equation}
			\end{lemma}
			
			\vspace*{1mm}
			\begin{proof}
				Refer to Lemma~4.5, 4.6 of \cite{Babuska-Nobile-Tempone}.
			\end{proof}

				\vspace*{1mm}
		\begin{remark}[analytic extension and operator independence]
		\label{rem:analytic-solution}
			We note that if the solution of a general-form SBVP is analytically extendable in the stochastic variables, then the stochastic collocation scheme and its error analysis remain independent of the operator acting in the physical domain.
		\end{remark}
			
		\vspace*{1mm}
			With the essential tools in place, we are now prepared to present an a priori error estimate, which also serves as a convergence result, for the discrete solution $(u_{p,\mathbb{U}}, \boldsymbol{g}_{p,\mathbb{G}})$, as described in (\ref{def:discrete-solution}). This result is formalized in the following theorem.
			
			\vspace*{1mm}
			\begin{theorem}[a poriori error bound]
				\label{The:a-poriori-error-bound}
				Assume that the strong solution $u$ of (\ref{eq:nondivergence-random-variable}) satisfies $u(\boldsymbol{y}, \cdot) \in H^{\varrho+2}(\mathcal{D})$ for a.e. $\boldsymbol{y} \in \boldsymbol{\Gamma}$, with some real $\varrho > 0$. Under Assumption~\ref{assum:convergence-A-f}, there exist constants $C_{\ref{ineq:a-poriori-error-bound}}>0$ and $r_n>0$ for $n = 1, \cdots, N$, independent of the discretization parameters  $h$, $k$, and $p$, such that
				\begin{equation}
					\label{ineq:a-poriori-error-bound}
					\left\|  (u, \boldsymbol{g}) - (u_{\boldsymbol{p},\mathbb{U}}, \boldsymbol{g}_{\boldsymbol{p},\mathbb{G}})\right\| _{ L^2_\rho(\boldsymbol{\Gamma}) \otimes H^1_h(\mathcal{D}) }
					\leq 
					C_{\ref{ineq:a-poriori-error-bound}}
					\left(
					h^{min\lbrace k, \varrho \rbrace } 
					+
					\sum_{n=1}^{N}{\beta_n \exp\left( -r_n p_n^{\theta_n} \right) }
					\right), 
				\end{equation}
				where 
				\begin{itemize}
					\item 
					if $\Gamma_n$ is bounded \quad
					$\begin{cases}
						\theta_n = \beta_n = 1
						\\
						r_n = \log \left[  \frac{2\tau_n}{\lvert \Gamma_n \rvert} \left( 1+ \sqrt{1+ \frac{\lvert \Gamma_n \rvert^2}{4\tau_n^2}} \right) \right] 
					\end{cases}$
					\item 
					if $\Gamma_n$ is unbounded  \quad 
					$\begin{cases}
						\theta_n = \dfrac{1}{2}, ~~ \beta_n = O(\sqrt{p_n})
						\\
						r_n = \tau_n\delta_n,
					\end{cases}$
					\\
					where 
					%$\tau_n \in (0, \dfrac{1}{2\eta_n})$ 
					$\tau_n \in (0, {1}/{(2\eta_n)})$
					and $\delta_n$ is as considered in (\ref{eq:density-func-assum}) . 
				\end{itemize}
			\end{theorem}
			
			\vspace*{1mm}
			\begin{proof}
				Let$(\mathcal{I}_{\boldsymbol{p}}u, \mathcal{I}_{\boldsymbol{p}}\boldsymbol{g})$
				denote the stochastic discretization alone; indeed $(\mathcal{I}_{\boldsymbol{p}}u, \mathcal{I}_{\boldsymbol{p}}\boldsymbol{g})$ represents the interpolation of $(u, \boldsymbol{g})$ as described in (\ref{def:interpolation-op}). The triangle inequality implies that
				\begin{multline}
					\label{ineq:a-poriori-triangle}
					\left\|  (u, \boldsymbol{g}) - (u_{\boldsymbol{p},\mathbb{U}}, \boldsymbol{g}_{\boldsymbol{p},\mathbb{G}})\right\| _{ L^2_\rho(\boldsymbol{\Gamma}) \otimes H^1_h(\mathcal{D}) } 
					\\
					\leq 
					\left\|  (u, \boldsymbol{g}) - (\mathcal{I}_{\boldsymbol{p}} u, \mathcal{I}_{\boldsymbol{p}} \boldsymbol{g})\right\| _{ L^2_\rho(\boldsymbol{\Gamma}) \otimes H^1_h(\mathcal{D}) }
					+
					\left\|  (\mathcal{I}_{\boldsymbol{p}}u, \mathcal{I}_{\boldsymbol{p}}\boldsymbol{g}) - (u_{\boldsymbol{p},\mathbb{U}}, \boldsymbol{g}_{\boldsymbol{p},\mathbb{G}})\right\| _{L^2_\rho(\boldsymbol{\Gamma}) \otimes H^1_h(\mathcal{D})}.
				\end{multline}
				For the second term of the right-hand side, %of inequality (\ref{ineq:a-poriori-triangle}), 
				we recall that  $u_{\boldsymbol{p},\mathbb{U}} = \left( \mathcal{I}_{\boldsymbol{p}} u\right) _{\mathbb{U}}$,  $\boldsymbol{g}_{\boldsymbol{p} ,\mathbb{G}} = \left( \mathcal{I}_{\boldsymbol{p}} \boldsymbol{g}\right) _{\mathbb{G}}$ and moreover, that $\mathcal{I}_{\boldsymbol{p}} u$ has the same regularity as the exact solution $u$ with respect to the spatial variable $\boldsymbol{x}$. Hence, Theorem~\ref{lem:a-priori-error-estimate} implies that there exist a constant $C_{\ref{ineq:a-priori-error-estimate}}>0$ such that 
				\begin{equation}
					\label{ineq:a-poriori-error-space-interpolation}
					\left\|  (\mathcal{I}_{\boldsymbol{p}}u, \mathcal{I}_{\boldsymbol{p}} \boldsymbol{g}) - (u_{\boldsymbol{p} ,\mathbb{U}}, \boldsymbol{g}_{\boldsymbol{p},\mathbb{G}})\right\| _{L^2_\rho(\boldsymbol{\Gamma}) \otimes H^1_h(\mathcal{D})}
					\leq
					C_{\ref{ineq:a-priori-error-estimate}} h^{min\lbrace k, \varrho \rbrace } \left\| \mathcal{I}_{\boldsymbol{p}}u  \right\|_{ L^2_\rho(\boldsymbol{\Gamma}) \otimes H^{2+\varrho}(\mathcal{D}) }.
				\end{equation}
				From Lemma~\ref{lem:continuity-interpolation-op}, we obtain  
				\begin{equation}
					\label{ineq:a-poriori-error-space}
					\left\|  (\mathcal{I}_{\boldsymbol{p}} u, \mathcal{I}_{\boldsymbol{p}}\boldsymbol{g}) - (u_{\boldsymbol{p} ,\mathbb{U}}, \boldsymbol{g}_{\boldsymbol{p},\mathbb{G}})\right\| _{ L^2_\rho(\boldsymbol{\Gamma}) \otimes  H^1_h(\mathcal{D}) }
					\leq
					C_{\ref{ineq:a-priori-error-estimate}}
					C_{\ref{ineq:continuity-interpolation-op}} 
					\left\| u  \right\|_{C^0_{\varsigma}\left( \boldsymbol{\Gamma};H^{2+\varrho}(\mathcal{D})\right) } h^{min\lbrace k, \varrho \rbrace }.
				\end{equation}
				The first term on the right hand side of (\ref{ineq:a-poriori-triangle}) is an interpolation error where the synonyms $\boldsymbol{g} = \nabla u$, $\mathcal{I}_{\boldsymbol{p}} \boldsymbol{g} = \mathcal{I}_{\boldsymbol{p}} \nabla u$ lead to
				\begin{equation*}
					\left\|  (u, \boldsymbol{g}) - (\mathcal{I}_{\boldsymbol{p}} u, \mathcal{I}_{\boldsymbol{p}} \boldsymbol{g})\right\| _{ L^2_\rho(\boldsymbol{\Gamma}) \otimes H^1_h(\mathcal{D})}
					\sim
					\left\| u - \mathcal{I}_{\boldsymbol{p}} u \right\|_{L_\rho^2(\boldsymbol{\Gamma}) \otimes H^2(\mathcal{D})}. 
				\end{equation*} 
				To analyze this term, we use the inclusion (\ref{ineq:countiuity-embedding}), which ensures the existence of a constant $C>0$ such that
				%we start by passing form the norm $L_\rho^2$ to $L^2_{\hat{\rho}}$:
				\begin{equation*}
					\left\|  u - \mathcal{I}_{\boldsymbol{p}} u \right\|_{ L^2_\rho(\boldsymbol{\Gamma}) \otimes H^2(\mathcal{D})}
					\leq
					%\left\| \dfrac{\rho}{\hat{\rho}} \right\|_{L^\infty(\boldsymbol{\Gamma})}^{\frac{1}{2}}
					C
					\left\|  u - \mathcal{I}_{\boldsymbol{p}} u \right\| _{L^2_{\hat{\rho}}(\boldsymbol{\Gamma}) \otimes H^2(\mathcal{D})}.
				\end{equation*}
				As before, we indicate $\bullet_n$ as a quantity related to direction $y_n$ and $\bullet_n^*$ the analogous quantity relative to all other directions $y_j, ~ j\neq n$. First, we are going to employ a one dimensional argument and focus on the direction $y_1$. We define an interpolation operator $\mathcal{I}_1: C_{\sigma1}^0\left( \Gamma_1; L_{\hat{\rho}_1^*}^2(\boldsymbol{\Gamma}_1^*) \otimes  H^2 \right)  \rightarrow L_{\hat{\rho}_1}^2 \left( \Gamma_1; L_{\hat{\rho}_1^*}^2(\boldsymbol{\Gamma}_1^*) \otimes H^2 \right)$ as
				\begin{equation*}
					\mathcal{I}_1 v\left( y_1, \boldsymbol{y}_1^*,\boldsymbol{x} \right)
					=
					\sum_{j=1}^{p_1+1}  
					v\left( y_{1,j}, \boldsymbol{y}_1^*,\boldsymbol{x} \right) 
					\ell_{1,j}(y_1).
				\end{equation*}
				We consider $u$ as a function of $y_1$ with values in the Hilbert space $V:= L_{\hat{\rho}_1^*}^2(\boldsymbol{\Gamma}_1^*) \otimes  H^2$; namely $u \in L_{\hat{\rho}_1}^2\left(\Gamma_1; V \right)$. We aim to provide an upper bound for $\left\| u - \mathcal{I}_1 u \right\|_{L^2_{\hat{\rho}_1}(\Gamma_1; V)} $.
				From Remark~\ref{rem:regularity-w.r-paramter} we know that $u \in C_{\sigma_1}^0(\Gamma_1; V)$ and then, Lemma~\ref{lem:optimality-interpolation} implies that 
				\begin{equation*}
					\left\| u - \mathcal{I}_1 u \right\|_{L^2_{\hat{\rho}_1}(\Gamma_1; V)} 
					\leq
					C_{\ref{ineq:optimality-interpolation}}
					%{\color{red}
					\inf_{w \in \mathcal{P}_{p_1}(\Gamma_1) \otimes V}
					\left\| u - w \right\|_{C_{\varsigma_1}^0(\Gamma_1, V)}
					%}
					.
				\end{equation*}
				In the case of bounded $\Gamma_1$ we employ Lemma~\ref{lem:bonded-domain-min-approx} and in the case of unbounded $\Gamma_1$, we employ Lemma~\ref{lem:unbonded-domain-min-approx}, for the analytic solution $u$ comes from Lemma~\ref{lem:analytic-exte}. Hence, we obtain
				\begin{equation}
					\label{ineq:upper-bound-interpolation-1D}
					\left\| u - \mathcal{I}_1 u \right\|_{L_{\hat{\rho}_1}^2(\Gamma_1; V)} 
					\leq
					\begin{cases}
						Ce^{-r_1p_1} & \text{if $\Gamma_1$ is bounded},
						\\
						C_{\ref{ineq:unbonded-domain-min-approx}}O\left(\sqrt{p_1} \right) e^{-r_1\sqrt{p_1}} 
						& \text{if $\Gamma_1$ is unbounded},
					\end{cases}	
				\end{equation}
				where $r_1$, $C$ and $C_{\ref{ineq:unbonded-domain-min-approx}}$ being specified from Lemma~\ref{lem:bonded-domain-min-approx} and Lemma~\ref{lem:unbonded-domain-min-approx}.
				Now, let us return to the general interpolation operator $\mathcal{I}_{\boldsymbol{p}}$, which is considered as the composition of two interpolants: $\mathcal{I}_{\boldsymbol{p}}= \mathcal{I}_1 \circ \mathcal{I}_1^*$, where $\mathcal{I}_1^*$ is the interpolation operator in all directions $y_2,y_3,\dots,y_N$ except $y_1$; namely $\mathcal{I}_1^*: C_{\sigma_{1}^*}^0\left( \boldsymbol{\Gamma}_1^* ;H^2 \right) \rightarrow L_{\hat{\rho}_1^*}^2 \left( \boldsymbol{\Gamma}_1^*; H^2 \right)$. Thus, we have
				\begin{equation*}
					\label{ineq:error-interpolation-triangle}
					\left\| u - \mathcal{I}_{\boldsymbol{p}} u \right\|_{L_{\hat{\rho}}^2(\boldsymbol{\Gamma}) \otimes H^2(\mathcal{D})}
					\leq
					\left\| u - \mathcal{I}_1u \right\|_{L_{\hat{\rho}}^2(\boldsymbol{\Gamma}) \otimes H^2(\mathcal{D})}
					+ 
					\left\| \mathcal{I}_1 \left( u - \mathcal{I}_1^* u \right) 
					\right\|_{L_{\hat{\rho}}^2(\boldsymbol{\Gamma}) \otimes H^2(\mathcal{D})} .
				\end{equation*}
				The first term on the right-hand side is bounded by (\ref{ineq:upper-bound-interpolation-1D}). To bound the second term, from Lemma~\ref{lem:continuity-interpolation-op}, we have
				\begin{equation*}
					\label{ineq:continuity-interpolation1}
					\left\| \mathcal{I}_1 \left( u - \mathcal{I}_1^* u \right) 
					\right\|_{L_{\hat{\rho}}^2(\boldsymbol{\Gamma}) \otimes H^2(\mathcal{D})}
					\leq
					C_{\ref{ineq:continuity-interpolation-op}}
					\left\|  u - \mathcal{I}_1^* u \right\|_{C_{\varsigma_1}^0(\Gamma_1; V)}.
				\end{equation*}
				This is again an interpolation error. To bound this error, we proceed iteratively by defining an interpolation operator $\mathcal{I}_2$ for the next direction $y_2$, and then bound the resulting error in the direction $y_2$, and so on for all the remaining directions $y_3, \cdots, y_N$. Thus, we can iteratively bound the interpolation error in each direction to ultimately estimate the total interpolation error across all dimensions.
			\end{proof}
		
		\vspace*{2mm}
		\section{Numerical Experiment}
		\label{sec:Numerical-Experiment}
		
\begin{figure}[tbhp]
			\centering
			\subfloat[\scriptsize Gaussian random variable $\boldsymbol{y}$ and ${\boldsymbol{p}}=(p_1,8)$]{\label{fig:test1-error-exp-Gaussian-p1}\includegraphics[scale=.47]{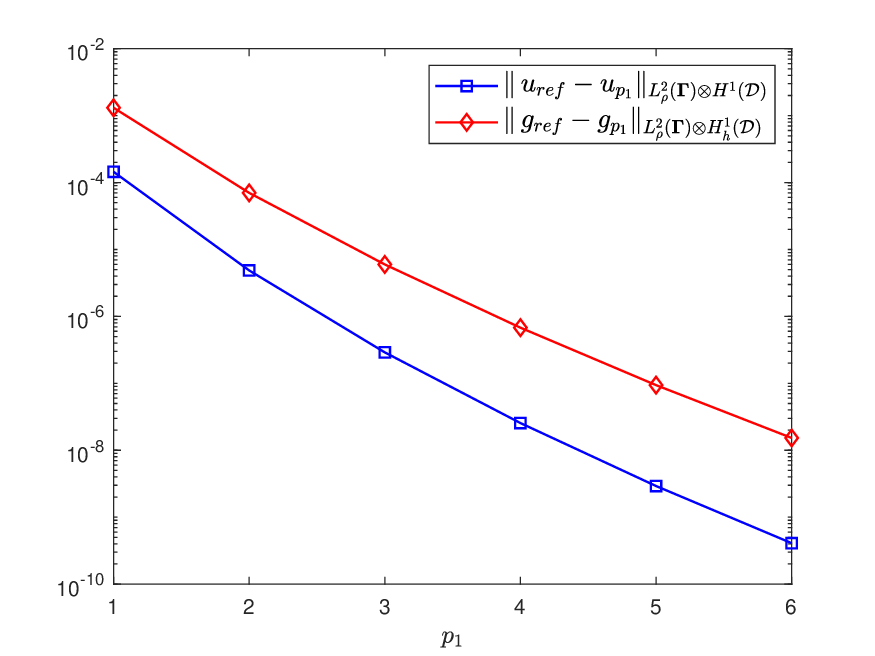}}
			\subfloat[\scriptsize Uniform random variable $\boldsymbol{y}$ and ${\boldsymbol{p}}=(p_1, 8)$]{\label{fig:test1-error-exp-Uniform-p1}\includegraphics[scale=.47]{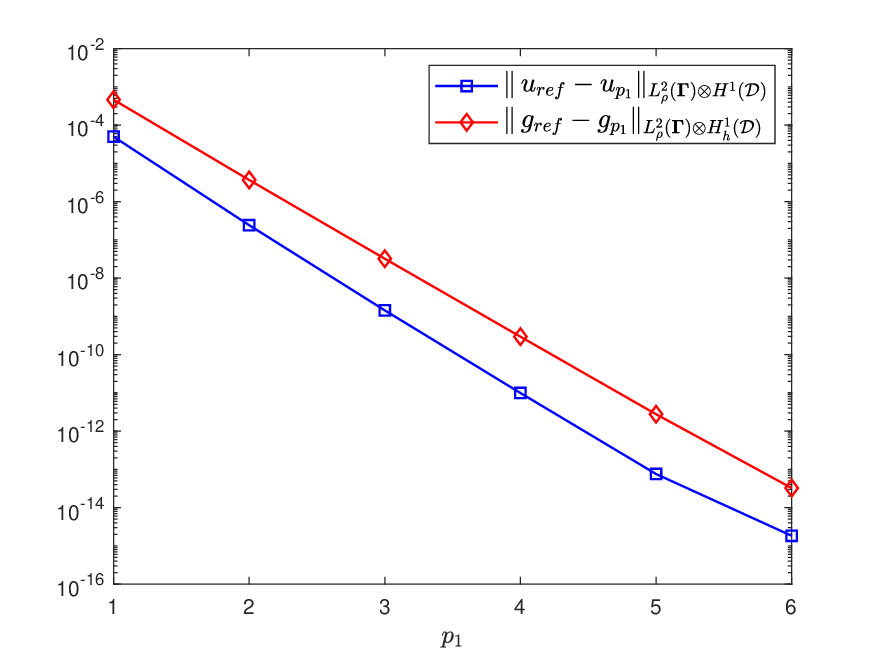}}	
			%\label{fig:test1-error-exp-Uniform-p1}
			\\
			\subfloat[\scriptsize Gaussian random variable $\boldsymbol{y}$ and ${\boldsymbol{p}}=(8,p_2)$]{\label{fig:test1-error-exp-Gaussian-p2}\includegraphics[scale=.47]{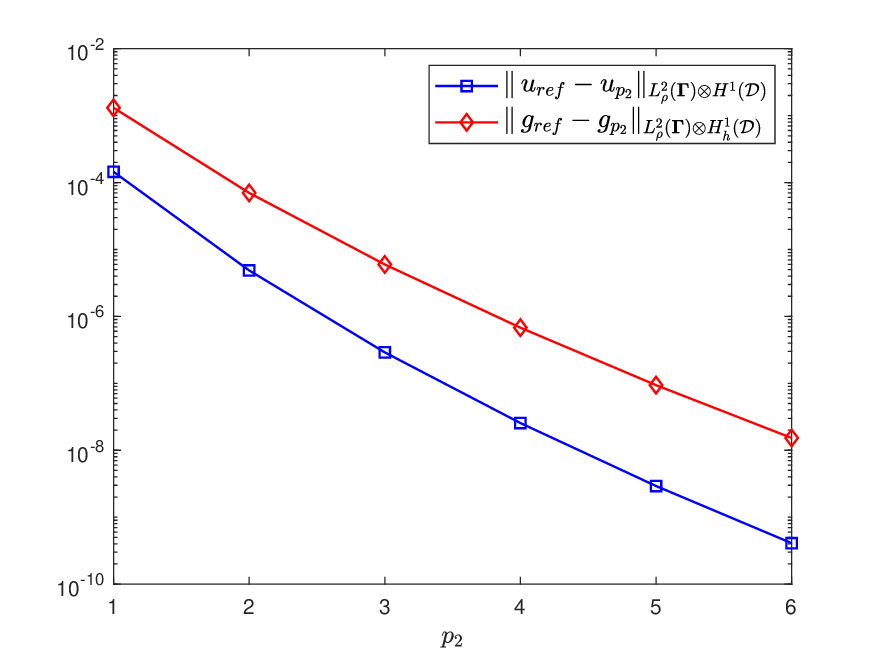}}
			\subfloat[\scriptsize Uniform random variable $\boldsymbol{y}$ and ${\boldsymbol{p}}=(8,p_2)$]{\label{fig:test1-error-exp-Uniform-p2}\includegraphics[scale=.47]{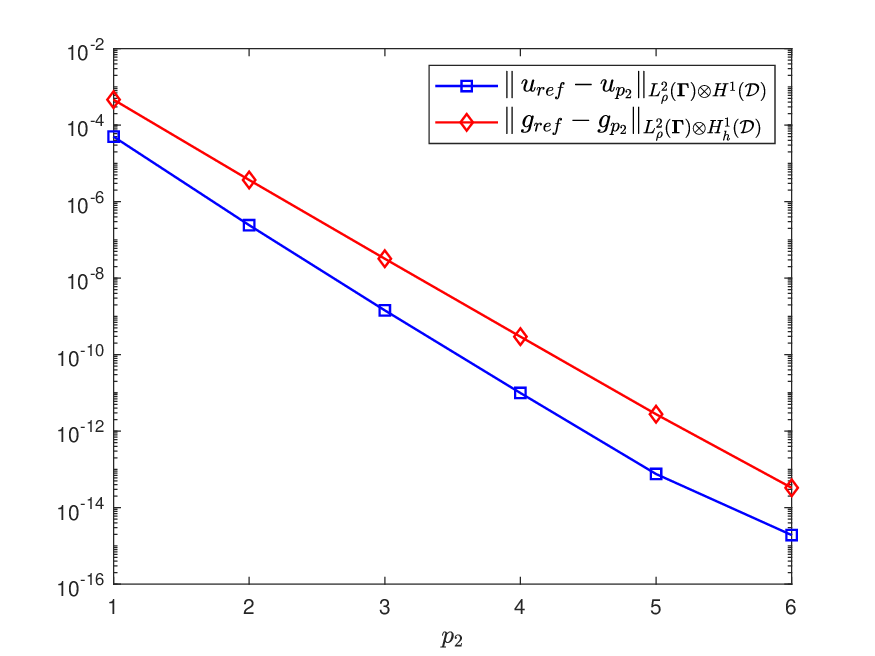}}
			\caption{Logarithm of the error versus the polynomial degree $p_n$, of the random variable $y_n$, for $n=1,2$.}
			\label{fig:test1-error-exp}
		\end{figure}

\begin{figure}[tbhp]
	\centering
	\subfloat[\scriptsize Gaussian random variable $\boldsymbol{y}$ and $\mathcal{P}_1$-element]{\label{fig:test1-error-poly-Gaussian-linear}\includegraphics[scale=.47]{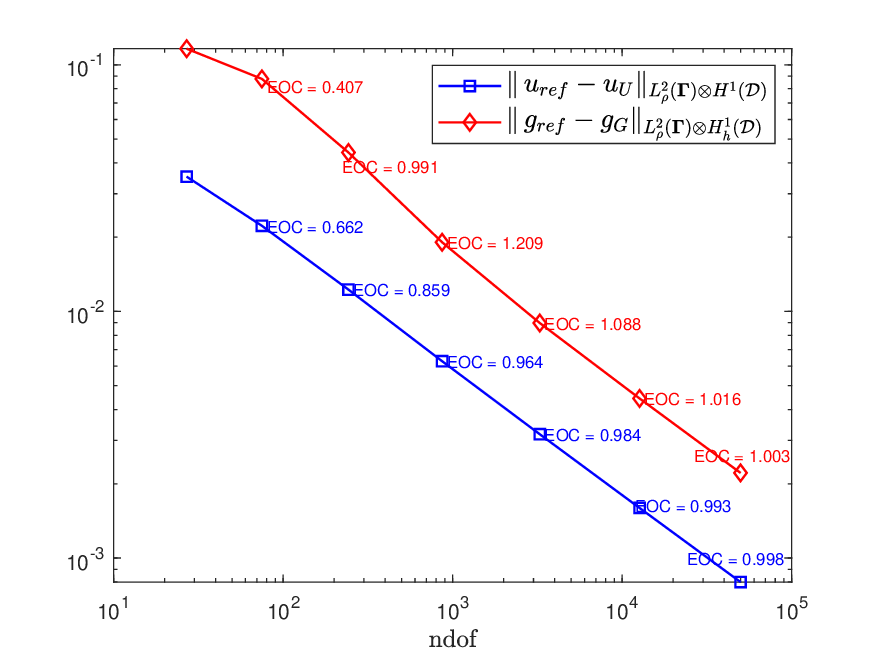}}
	\subfloat[\scriptsize Uniform random variable $\boldsymbol{y}$ and $\mathcal{P}_1$-element]{\label{fig:test1-error-poly-Uniform-linear}\includegraphics[scale=.47]{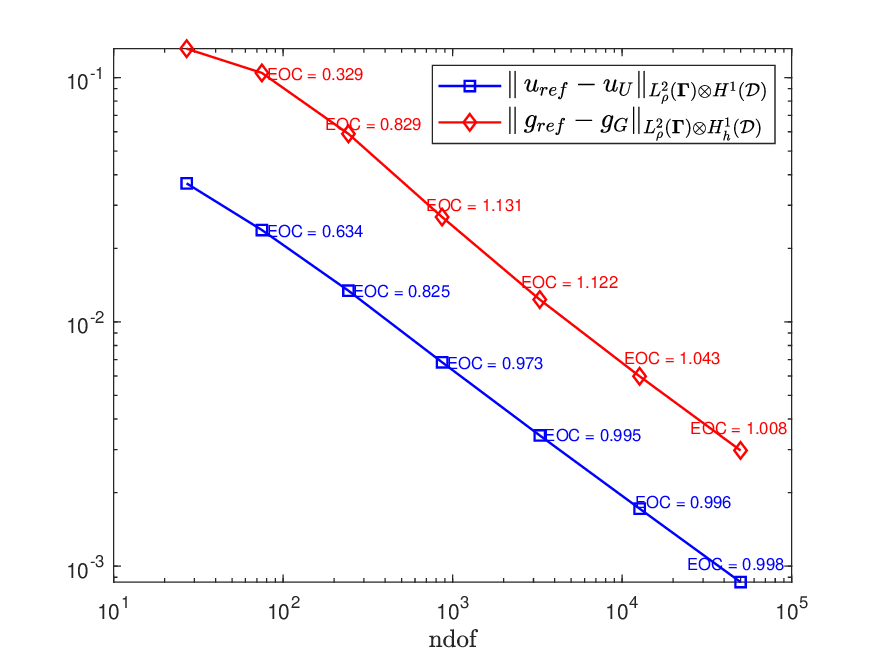}}	
	%\label{fig:test1-error-exp-Uniform-p1}
	\\
	\subfloat[{\scriptsize Gaussian random variable $\boldsymbol{y}$ and $\mathcal{P}_2$-element}]{\label{fig:test1-error-poly-Gaussian-quad}\includegraphics[scale=.47]{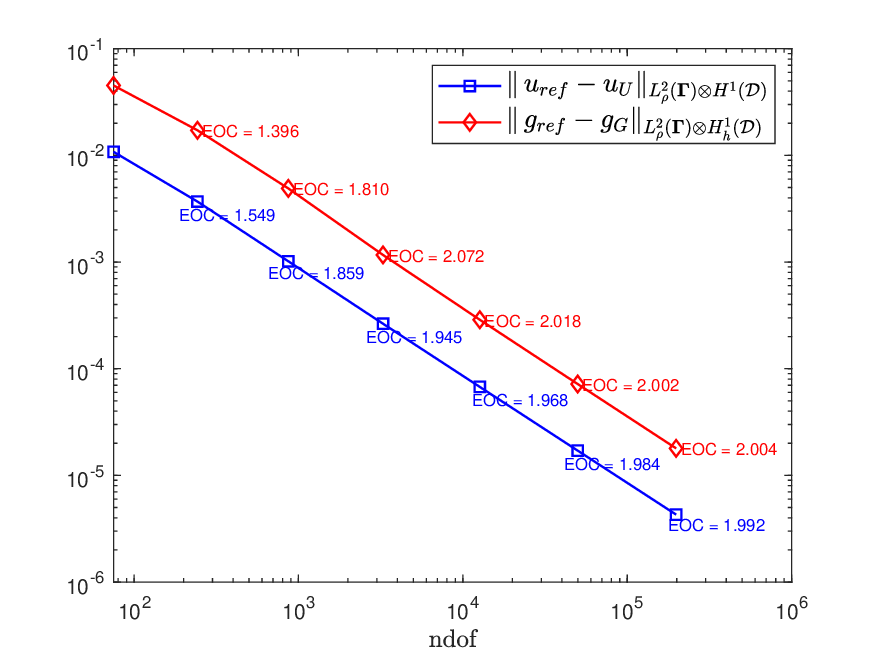}}
	\subfloat[{\scriptsize Uniform random variable $\boldsymbol{y}$ and $\mathcal{P}_2$-element}]{\label{fig:test1-error-poly-Uniform-quad}\includegraphics[scale=.47]{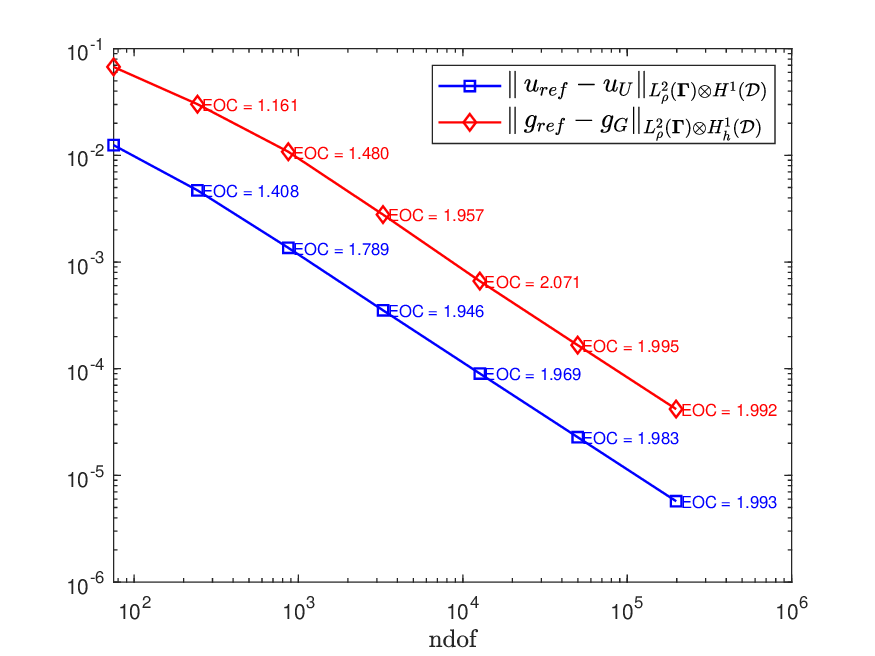}}
	\caption{(log–log plot) Error versus number of degrees of freedom (ndof), showing the convergence rates for $\mathcal{P}_k$-elements with $k=1,2$.}
	\label{fig:test1-error-polynomial}
\end{figure}

		%\subsection{Test problem with smooth data}
		%\label{test-problem-smooth-data}
		In this test problem, let $\mathcal{D} = (-1, 1) \times (-1,1)$. %and $f={(2-x_1^2 -x_2^2)}/{8}$. 
		For $\boldsymbol{x}=(x_1,x_2)$ and a random vector $\boldsymbol{y}=(y_1,y_2)$, where $y_1$ and $y_2$ are independent and identically distributed (i.i.d), the random diffusion coefficient matrix and the right-hand side are defined as
		\begin{align*}	
			\boldsymbol{A}(\boldsymbol{x}, \boldsymbol{y}) =&
			\begin{bmatrix}
				\frac{5}{2} + e^{ -\frac{y_1^2}{100}}  \left( \cos(\pi x_1) + \sin(\pi x_2) \right) 
				& 0
				\\
				0 &
				\frac{5}{2} + e^{-\frac{y_2^2}{100} }  \left( \sin(\pi x_1) + \cos(\pi x_2) \right)
			\end{bmatrix},
			\\
			f(\boldsymbol{x}, \boldsymbol{y}) =& \frac{1}{8}(2-x_1^2 -x_2^2).
		\end{align*}
		Since the exact solution of (\ref{eq:nondivergence-random-variable}) is unknown for the given data, the error is computed using a reference solution. The pair $(u_{\text{ref}}, \boldsymbol{g}_{\text{ref}})$ represents a reference solution obtained by solving (\ref{def:discrete-solution}) with ${\boldsymbol{p}}=(8,8)$ ($N_{\boldsymbol{p}}=81$) and polynomial degree of $k =4$, with a mesh size of $h = {1}/{128}$.
		We test two different types of distributions for the random variable $\boldsymbol{y}$: Gaussian and Uniform, which have unbounded and bounded domains, respectively. The collocation points corresponding to these distributions are determined by the Cartesian products of the roots of Hermite polynomials for the Gaussian case, and Legendre polynomials for the Uniform case. 
		To demonstrate the performance of the discretization scheme proposed in Section~\ref{sec:discretization}, we test its effectiveness separately in both the stochastic and physical domains. In Figure~\ref{fig:test1-error-exp}, we present the logarithm of the discretization error in the stochastic domain, while keeping the polynomial degree of the finite element space defined in (\ref{def:Galerkin-spaces}) fixed at $k=4$, and the mesh size in the physical domain fixed at $h={1}/{128}$. The focus here is on varying the number of collocation points and, consequently, the degree of the Lagrange polynomials with respect to one random variable, while keeping the number of collocation points fixed, equal to $9$, for the other random variable. This approach aims to assess how the approximation accuracy changes as the number of collocation points for one random variable increases. We observe an exponential rate of convergence with respect to the polynomial degree used for the approximation. In Figure~\ref{fig:test1-error-polynomial}, we report the (log–log) error of the discretization in the physical domain, keeping the number of collocation points for both random variables fixed at ${\boldsymbol{p}}=(8,8)$, for both Gaussian and Uniform distributions.
%{Logarithm of the $H^1$-error of the mean value  vs. the polynomial degree $p_n$, of the random variable $y_n$, for $n=1,2$.}
%	

%
%\begin{figure}[tbhp]
%	\centering
%	\begin{subfigure}[b]{0.45\textwidth}%\centering
%		%\rule{2cm}{2cm}
%		\includegraphics[scale=.55]{Figures/test1-Gaussian-exp-p1-s}
%		\caption{Gaussian random variable $y_1$}
%		\label{fig:test1-error-exp-Gaussian-p1}
%	\end{subfigure}
%	\begin{subfigure}[b]{0.45\textwidth}%\centering
%		\includegraphics[scale=.55]{Figures/test1-Uniform-exp-p1-s}
%		\caption{Uniform random variable $y_1$}
%		\label{fig:test1-error-exp-Uniform-p1}
%	\end{subfigure}
%	\\*[4mm]
%	\begin{subfigure}[b]{0.45\textwidth}%\centering
%		\includegraphics[scale=.55]{Figures/test1-Gaussian-exp-p2-s}
%		\caption{Gaussian random variable $y_2$}
%		\label{fig:test1-error-exp-Gaussian-p2}
%	\end{subfigure}
%	\begin{subfigure}[b]{0.45\textwidth}%\centering
%		\includegraphics[scale=.55]{Figures/test1-Uniform-exp-p2-s}
%		\caption{Uniform random variable $y_2$}
%		\label{fig:test1-error-exp-Uniform-p2}
%	\end{subfigure}
%	\caption{Logarithm of the $H^1$-error of the mean value  vs. the polynomial degree $p_n$, of the random variable $y_n$, for $n=1,2$.}
%	\label{fig:test1-error-exp}
%\end{figure}

%%%%%%%%%%%%%%%%%%%%%%%%%%%%%%%%%%%%%%%%%%%%%%%%%%%%%%%%%%%%%%%%%%%%%%%%%%%%%%%%%%%
%%%%%%%%%%%%%%%%%%%%%%%%%%%%%%%%%%%%%%%%%%%%%%%%%%%%%%%%%%%%%%%%%%%%%%%%%%%%%%%%%%%

	\bibliographystyle{siamplain}
	\bibliography{references}
	%
	%\bigskip
	%\medskip
	
	\vspace*{3mm}
	\hspace*{.8mm}
	\noindent {\footnotesize Amireh Mousavi, 
	Friedrich-Schiller-Universität Jena, 07743, Jena, Germany 
	\vspace*{-0.5mm}
	\\
	\hspace*{7mm}
	\textit{Email address:} \texttt{amireh.mousavi@gmail.com}}
	
\end{document}